\documentclass{amsart}
 \usepackage[sorted-cites]{amsrefs}
\usepackage{amsfonts}
\usepackage{graphicx}
\usepackage{amscd}
\usepackage{amsmath}
\usepackage{amssymb}
\usepackage{tikz}
\usepackage{pgfplots,pgfplotstable}
\usetikzlibrary{arrows.meta,fit}
\usetikzlibrary{arrows.meta,calc}
\definecolor{base}{RGB}{74,111,165}
\definecolor{lev}{RGB}{127,156,200}
\definecolor{fib}{RGB}{217,130,43}
\definecolor{sphfill}{RGB}{223,231,242}
\usepackage{xcolor}
\usepackage{hyperref}
\usepackage{enumitem}
\newtheorem{theorem}{\bf Theorem}[section]
\newtheorem{lemma}{\bf Lemma}[section]
\newtheorem{proposition}{\bf Proposition}[section]
\newtheorem{corollary}{\bf Corollary}[section]

\newtheorem{remark}{\bf Remark}
\makeatletter \@addtoreset{equation}{section} \makeatother

\def \k{\tilde{g}}

\newcommand{\Ric}{\operatorname{Ric}}
\newcommand{\Hess}{\nabla^2}
\newcommand{\Crit}{\operatorname{Crit}}
\newcommand{\R}{\mathbb R}
\newcommand{\C}{\mathbb C}
\newcommand{\CP}{\mathbb{CP}}
\newcommand{\OO}{\mathcal O}

\begin{document}

\title[Steady conformally K\"{a}hler Ricci solitons]{Steady conformally K\"{a}hler Ricci solitons\\ with isoparametric potential}

\author{Ernani Ribeiro Jr}
\author{Hung Tran}

\address{Departamento  de Matem\'atica, Universidade Federal do Cear\'a, 60455-760, Fortaleza, Brazil}
\email{ernani@mat.ufc.br}

\address{Department of Mathematics and Statistics, Texas Tech University, Lubbock, TX 79409} \email{hung.tran@ttu.edu}

\thanks{E. Ribeiro Jr was partially supported by CNPq/Brazil - 305128/2025-6 and 351492/2025-9.}

\keywords{Steady Ricci soliton; Ricci flow; conformally K\"ahler metric; isoparametric potential} \subjclass[2020]{Primary 53C25; Secondary 53E20, 53C21, 53C55}

\date{\today}

\begin{abstract}
We classify complete, non-compact, simply connected, irreducible
four-dimensional steady gradient Ricci solitons that are globally
conformally K\"ahler and have an isoparametric potential. The steady Ricci solitons are $U(2)$-invariant nuts on
$\mathbb R^4$, including the Bryant and Cao solitons, and bolts on
$\mathbb{S}^2 \times \mathbb R^2$ and on $\mathcal O(-m)\to\mathbb{CP}^1$. For
each $m\ge3$ there is exactly one non-collapsed bolt. A key step, valid for Ricci solitons of every sign, is that away from the
locally conformally flat case, the conformal factor is aligned with the
potential function, which yields cohomogeneity-one
symmetry. 
\end{abstract}

\maketitle

\section{Introduction}\label{intro}

A Riemannian manifold $(M^n,\,g)$ is called a \textit{gradient Ricci soliton} (GRS) if there exists a smooth potential function $f$ such that the Ricci tensor $\Ric$ satisfies 
\begin{equation}\label{grs}
 \Ric+\Hess f=\lambda g
\end{equation} for some constant $\lambda\in \mathbb{R}.$ Here $\Hess\, f$ denotes the Hessian of $f$. It is called expanding if $\lambda<0,$ steady if $\lambda=0,$ or shrinking if $\lambda>0.$ GRS generate self-similar solutions to the Ricci flow, introduced by
R. Hamilton~\cite{hamilton82three}, and play a central role in the analysis of its
singularities, as confirmed in every dimension by Bamler's structure theory for
non-collapsed limits of Ricci flows~\cites{bamler23compactness,bamler21structure}.
The theory of GRS has several celebrated applications, including
\cites{bohm08manifolds, brendle08classification, brendle09manifolds, cao23curvature,
hamilton86four, hamilton95formation, perelman02entropy, perelman03finite,
perelman03surgery}; see \cites{CaoSurvey, cao24survey} for surveys.

Steady GRS arise naturally as Type~II singularity models of the Ricci flow~\cite{cao97limits,hamilton95formation}, and their classification is a fundamental problem. In dimension two, Hamilton~\cite{hamilton88surfaces}
showed that every complete steady GRS is either flat or isometric, up to scaling, to the cigar soliton. In dimension three, Brendle~\cite{brendle13rotational} proved that the Bryant soliton is the only non-flat, non-collapsed
example, and Lai~\cite{Lai2025FlyingWings, Lai2025} constructed the collapsed flying wings and established their $O(2)$-symmetry. In
dimensions $n\ge4$, the problem remains largely open. In dimension four, Bamler, Chow, Deng, Ma and Zhang~\cite{BCDMZ22} studied steady
soliton singularity models with $3$-cylindrical tangent flows at infinity. Most higher-dimensional rigidity results rely on suitable curvature assumptions---including local conformal flatness~\cite{CC, CM}, harmonic or half-harmonic Weyl curvature~\cite{CaoXie2026, kim17classification, Li2025Rigidity}, vanishing Bach or $D$-tensors~\cite{cao2014bach, CaoYu}, or curvature decay~\cite{CH2016, deng2020higher, MMS2026}---which typically reduce the metric to rotational symmetry; see also \cite{MSW19,MW11} for related results.

Stepping down the symmetry hierarchy, a natural alternative is to impose a geometric structure on the potential function $f$. A smooth function is \emph{transnormal} if $|\nabla f|$ is constant on its level sets, and \emph{isoparametric} if its regular level sets additionally have constant mean curvature---a classical notion rooted in geometric optics~\cite{cartan38familles,levicivita37famiglie,somigliana18relazione} whose classification on round spheres (Yau's Question~34~\cite{yau93open}) has a rich history~\cite{abresch83isoparametric,cecil07isoparametric,chi20isoparametric,chi20story,munzner80isoparametrische,nomizu73some}. On a GRS, transnormality and isoparametricity are equivalent~\cite{Tran}, and the notion appears
throughout the theory. All known non-trivial examples constructed with cohomogeneity one symmetry have isoparametric potential~\cite{buttsworth21su2,Cao,cao97limits,chave96class, DW, feldman03rotationally,koiso90rotationally,pedersen99quasi}.
It is used in the classification of four-dimensional GRS with
harmonic Weyl curvature~\cite{kim17classification} or constant scalar curvature~\cite{cheng23rigidity,fernandez16constant}.
It also appears asymptotically along expanding K\"ahler GRS~\cite[Section~3]{conlon24classification}. In the K\"ahler setting, complete four-dimensional GRS with isoparametric potential are classified
in~\cite[Theorem~A]{Tran}; the simply connected ones turn out to have maximal symmetry as characterized in~\cite{tran23kahler}. It is also part of a dichotomy discussed in ~\cite{tran24toric}; the Riemannian case is addressed in \cite{TY}.


In dimension four, GRS display a strong rigidity toward the K\"ahler category; to our
knowledge, all known examples of nontrivial (irreducible and non-Einstein)
four-dimensional shrinking GRS are K\"ahler; see \cite{li26kahler, CRT26}. If
irreducibility is dropped, every known four-dimensional shrinking GRS is either
K\"ahler or has constant scalar curvature; see \cite{li26kahler, cheng23rigidity,
fernandez16constant, CRT26}. In the latter case, (\ref{eq3a}) shows that the potential
function is transnormal, hence isoparametric~\cite{Tran}. For steady GRS, constant
scalar curvature forces $\mathrm{Ric}\equiv0$ by (\ref{eq4a}), so the isoparametric
condition is its natural non-trivial relaxation, while the K\"ahler condition excludes the Bryant soliton~\cite{Bryant}, which is conformally flat and hence conformally K\"ahler. This raises the question of whether a rigidity phenomenon persists in the larger class
of \emph{conformally K\"ahler} metrics, which contains both the K\"ahler steady
solitons and the Bryant soliton, and thus offers a natural framework for bridging
K\"ahler and Riemannian GRS.


To state our results, we use the following leafwise transnormal hypothesis: $|\nabla f|^2$ is constant on each connected component of each regular level of the potential function. 
Furthermore, a metric is said to be conformally K\"ahler if there is a smooth function $u$ and a complex structure $J$ defined on $M$ such that
\begin{equation}\label{ck}
 e^{2u}g\quad\text{is K\"{a}hler with respect to }J.
\end{equation}

Our main result assumes no curvature or asymptotic conditions. 

\begin{theorem}\label{classification}
Let $(M^4,g,f)$ be a complete, connected, non-compact, simply connected, irreducible steady GRS with non-constant leafwise transnormal potential and whose metric is conformally K\"ahler.  Then, up to isometry, homothety, and addition of a constant to $f$, it belongs to exactly one of the following families:
\begin{center}
\renewcommand{\arraystretch}{1.15}
\setlength{\tabcolsep}{3pt}
\begin{tabular}{lll}
\hline
{\rm Endpoint} & {\rm Manifold} & {\rm Complete normalized parameters}\\
\hline
Nut & $\R^4$ & $\delta\in(-\infty,0]$ \quad ($\delta=0$: Bryant; $\delta=-1/2$: Cao)\\
Twisted bolt, $m=1,2$ & $\OO(-m)$ & $\mu\in(-\infty,0)$ \quad ($m=2$: K\"ahler)\\
Twisted bolt, $m\geq3$ & $\OO(-m)$ & $\mu\in(-\infty,\mu_m]$, with a unique $\mu_m<0$\\
Untwisted bolt & $\mathbb{S}^2 \times \mathbb R^2$ & $\mu\in(-\infty,0)$\\
\hline
\end{tabular}
\end{center}
Conversely, every listed trajectory satisfies all the hypotheses. 
\end{theorem}
The models above use the distance-coordinate metric
\begin{equation}\label{ansatz}
 g=dt^2+p(t)^2\eta^2+q(t)^2h,\qquad f=f(t),\qquad
 d\eta=c\omega_h,\qquad K_h=4,
\end{equation}
where $h$ is a metric on a surface with Gaussian curvature $K_h$ and area form $\omega_h$, and $\eta$ denotes a connection one-form on a circle bundle over that surface. 
\begin{itemize}
    \item A \emph{nut} is a smooth point at $t=0$, with $p'(0)=q'(0)=1$. 
    \item A \emph{bolt} is a smooth two-sphere at $t=0$, with $p(0)=0$ and $q(0)>0$.
\end{itemize}  The nut-bolt terminology is inherited from \cite{gibbons79classification}, as in \cite{FT26a, FT26b}.

 For a twisted model, $c=2$ and $\eta^2+h$ is locally the round unit three-sphere metric. In twisted degree $m\geq1$, the fiber period is $2\pi/m$ and $p'(0)=m$; the completion is $\OO(-m)\to\CP^1$, with orientation disregarded. The untwisted model has $c=0$, fiber period $2\pi$, and $p'(0)=1$. We normalize a nut by $f(0)=0$ and $f''(0)=-1$, a bolt by $q(0)=1$ and $f(0)=0$, and put
\[
 \delta=p^{(3)}(0)-q^{(3)}(0), \qquad \mu=f''(0).
\]

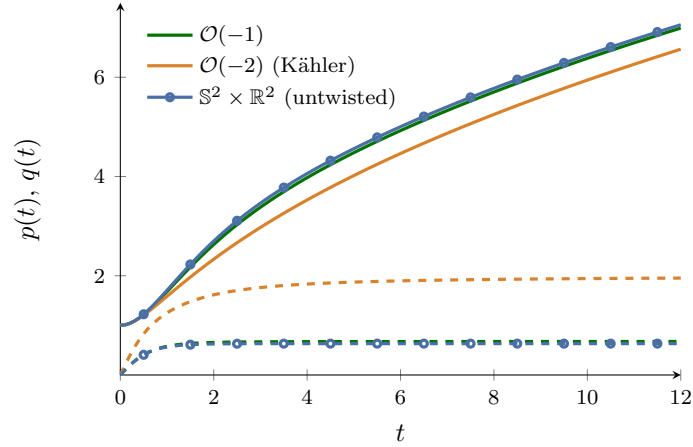
\begin{figure}[ht]
  \centering
\pgfplotstableread{
 t p q
0.0000 0.00000 1.00000
0.1000 0.09902 1.00994
0.2000 0.19246 1.03904
0.3000 0.27625 1.08540
0.4000 0.34837 1.14640
0.5000 0.40868 1.21923
0.6000 0.45815 1.30124
0.7000 0.49828 1.39012
0.8000 0.53062 1.48390
0.9000 0.55662 1.58098
1.0000 0.57753 1.68004
1.1000 0.59437 1.78002
1.2000 0.60796 1.88009
1.3000 0.61896 1.97959
1.4000 0.62789 2.07798
1.5000 0.63517 2.17488
1.6000 0.64113 2.27001
1.7000 0.64602 2.36314
1.8000 0.65005 2.45415
1.9000 0.65338 2.54295
2.0000 0.65615 2.62950
2.1000 0.65847 2.71381
2.2000 0.66041 2.79590
2.3000 0.66204 2.87581
2.4000 0.66342 2.95362
2.5000 0.66460 3.02938
2.6000 0.66560 3.10318
2.7000 0.66646 3.17511
2.8000 0.66721 3.24525
2.9000 0.66785 3.31369
3.0000 0.66842 3.38052
3.1000 0.66891 3.44581
3.2000 0.66935 3.50965
3.3000 0.66973 3.57211
3.4000 0.67008 3.63327
3.5000 0.67039 3.69320
3.6000 0.67067 3.75197
3.7000 0.67092 3.80963
3.8000 0.67115 3.86624
3.9000 0.67136 3.92185
4.0000 0.67155 3.97652
4.1000 0.67173 4.03030
4.2000 0.67189 4.08322
4.3000 0.67204 4.13533
4.4000 0.67219 4.18667
4.5000 0.67232 4.23726
4.6000 0.67244 4.28715
4.7000 0.67256 4.33637
4.8000 0.67267 4.38494
4.9000 0.67277 4.43289
5.0000 0.67287 4.48024
5.1000 0.67297 4.52703
5.2000 0.67306 4.57326
5.3000 0.67314 4.61896
5.4000 0.67322 4.66416
5.5000 0.67330 4.70886
5.6000 0.67337 4.75309
5.7000 0.67344 4.79686
5.8000 0.67351 4.84019
5.9000 0.67358 4.88309
6.0000 0.67364 4.92557
6.1000 0.67370 4.96765
6.2000 0.67376 5.00934
6.3000 0.67381 5.05064
6.4000 0.67387 5.09159
6.5000 0.67392 5.13217
6.6000 0.67397 5.17240
6.7000 0.67402 5.21230
6.8000 0.67406 5.25186
6.9000 0.67411 5.29110
7.0000 0.67415 5.33003
7.1000 0.67419 5.36866
7.2000 0.67424 5.40698
7.3000 0.67428 5.44502
7.4000 0.67431 5.48277
7.5000 0.67435 5.52024
7.6000 0.67439 5.55744
7.7000 0.67442 5.59438
7.8000 0.67446 5.63106
7.9000 0.67449 5.66748
8.0000 0.67453 5.70365
8.1000 0.67456 5.73959
8.2000 0.67459 5.77528
8.3000 0.67462 5.81074
8.4000 0.67465 5.84597
8.5000 0.67468 5.88098
8.6000 0.67471 5.91577
8.7000 0.67473 5.95034
8.8000 0.67476 5.98470
8.9000 0.67479 6.01885
9.0000 0.67481 6.05280
9.1000 0.67484 6.08655
9.2000 0.67486 6.12010
9.3000 0.67488 6.15346
9.4000 0.67491 6.18664
9.5000 0.67493 6.21962
9.6000 0.67495 6.25242
9.7000 0.67497 6.28504
9.8000 0.67500 6.31749
9.9000 0.67502 6.34976
10.0000 0.67504 6.38186
10.1000 0.67506 6.41379
10.2000 0.67508 6.44555
10.3000 0.67510 6.47715
10.4000 0.67511 6.50860
10.5000 0.67513 6.53988
10.6000 0.67515 6.57101
10.7000 0.67517 6.60198
10.8000 0.67519 6.63280
10.9000 0.67520 6.66348
11.0000 0.67522 6.69401
11.1000 0.67524 6.72439
11.2000 0.67525 6.75463
11.3000 0.67527 6.78473
11.4000 0.67528 6.81469
11.5000 0.67530 6.84452
11.6000 0.67531 6.87421
11.7000 0.67533 6.90376
11.8000 0.67534 6.93319
11.9000 0.67536 6.96249
12.0000 0.67537 6.99166
}\degone
\pgfplotstableread{
 t p q
0.0000 0.00000 1.00000
0.1000 0.19804 1.00990
0.2000 0.38517 1.03850
0.3000 0.55404 1.08293
0.4000 0.70191 1.13960
0.5000 0.82952 1.20505
0.6000 0.93926 1.27644
0.7000 1.03393 1.35161
0.8000 1.11612 1.42900
0.9000 1.18800 1.50752
1.0000 1.25130 1.58640
1.1000 1.30740 1.66513
1.2000 1.35740 1.74334
1.3000 1.40216 1.82079
1.4000 1.44240 1.89731
1.5000 1.47869 1.97281
1.6000 1.51153 2.04720
1.7000 1.54132 2.12046
1.8000 1.56840 2.19257
1.9000 1.59308 2.26353
2.0000 1.61562 2.33334
2.1000 1.63624 2.40202
2.2000 1.65515 2.46959
2.3000 1.67251 2.53607
2.4000 1.68848 2.60149
2.5000 1.70320 2.66589
2.6000 1.71680 2.72928
2.7000 1.72936 2.79170
2.8000 1.74101 2.85318
2.9000 1.75181 2.91375
3.0000 1.76185 2.97344
3.1000 1.77120 3.03227
3.2000 1.77992 3.09027
3.3000 1.78806 3.14747
3.4000 1.79567 3.20390
3.5000 1.80280 3.25957
3.6000 1.80950 3.31452
3.7000 1.81578 3.36877
3.8000 1.82169 3.42233
3.9000 1.82727 3.47523
4.0000 1.83252 3.52750
4.1000 1.83748 3.57914
4.2000 1.84218 3.63018
4.3000 1.84662 3.68064
4.4000 1.85084 3.73053
4.5000 1.85484 3.77987
4.6000 1.85864 3.82868
4.7000 1.86226 3.87696
4.8000 1.86570 3.92475
4.9000 1.86898 3.97204
5.0000 1.87211 4.01886
5.1000 1.87510 4.06521
5.2000 1.87796 4.11111
5.3000 1.88070 4.15658
5.4000 1.88332 4.20161
5.5000 1.88583 4.24623
5.6000 1.88824 4.29044
5.7000 1.89055 4.33425
5.8000 1.89277 4.37768
5.9000 1.89490 4.42073
6.0000 1.89696 4.46341
6.1000 1.89894 4.50573
6.2000 1.90084 4.54770
6.3000 1.90268 4.58933
6.4000 1.90445 4.63062
6.5000 1.90616 4.67159
6.6000 1.90781 4.71223
6.7000 1.90940 4.75256
6.8000 1.91094 4.79259
6.9000 1.91244 4.83231
7.0000 1.91388 4.87174
7.1000 1.91528 4.91088
7.2000 1.91663 4.94974
7.3000 1.91794 4.98833
7.4000 1.91921 5.02664
7.5000 1.92045 5.06469
7.6000 1.92165 5.10248
7.7000 1.92281 5.14002
7.8000 1.92394 5.17730
7.9000 1.92504 5.21434
8.0000 1.92610 5.25114
8.1000 1.92714 5.28770
8.2000 1.92815 5.32403
8.3000 1.92913 5.36013
8.4000 1.93009 5.39601
8.5000 1.93102 5.43167
8.6000 1.93193 5.46712
8.7000 1.93281 5.50235
8.8000 1.93367 5.53737
8.9000 1.93451 5.57219
9.0000 1.93533 5.60681
9.1000 1.93613 5.64123
9.2000 1.93691 5.67545
9.3000 1.93768 5.70948
9.4000 1.93842 5.74333
9.5000 1.93915 5.77699
9.6000 1.93986 5.81046
9.7000 1.94055 5.84376
9.8000 1.94123 5.87688
9.9000 1.94189 5.90982
10.0000 1.94254 5.94260
10.1000 1.94317 5.97520
10.2000 1.94380 6.00764
10.3000 1.94440 6.03991
10.4000 1.94500 6.07202
10.5000 1.94558 6.10398
10.6000 1.94615 6.13577
10.7000 1.94671 6.16741
10.8000 1.94726 6.19890
10.9000 1.94779 6.23024
11.0000 1.94832 6.26143
11.1000 1.94883 6.29247
11.2000 1.94934 6.32337
11.3000 1.94984 6.35413
11.4000 1.95032 6.38475
11.5000 1.95080 6.41522
11.6000 1.95127 6.44556
11.7000 1.95173 6.47577
11.8000 1.95218 6.50584
11.9000 1.95262 6.53578
12.0000 1.95305 6.56559
}\degtwo
\pgfplotstableread{
 t p q
0.0000 0.00000 1.00000
0.1000 0.09901 1.00995
0.2000 0.19242 1.03923
0.3000 0.27599 1.08622
0.4000 0.34753 1.14864
0.5000 0.40673 1.22385
0.6000 0.45456 1.30921
0.7000 0.49259 1.40227
0.8000 0.52254 1.50087
0.9000 0.54600 1.60318
1.0000 0.56434 1.70768
1.1000 0.57867 1.81315
1.2000 0.58987 1.91859
1.3000 0.59866 2.02323
1.4000 0.60555 2.12647
1.5000 0.61097 2.22789
1.6000 0.61525 2.32715
1.7000 0.61863 2.42405
1.8000 0.62130 2.51843
1.9000 0.62343 2.61024
2.0000 0.62511 2.69946
2.1000 0.62645 2.78609
2.2000 0.62753 2.87020
2.3000 0.62838 2.95186
2.4000 0.62906 3.03115
2.5000 0.62961 3.10817
2.6000 0.63005 3.18303
2.7000 0.63040 3.25583
2.8000 0.63068 3.32668
2.9000 0.63091 3.39569
3.0000 0.63109 3.46297
3.1000 0.63124 3.52860
3.2000 0.63136 3.59268
3.3000 0.63145 3.65531
3.4000 0.63153 3.71657
3.5000 0.63159 3.77654
3.6000 0.63164 3.83529
3.7000 0.63168 3.89288
3.8000 0.63171 3.94940
3.9000 0.63174 4.00488
4.0000 0.63176 4.05940
4.1000 0.63178 4.11300
4.2000 0.63179 4.16572
4.3000 0.63180 4.21762
4.4000 0.63181 4.26873
4.5000 0.63182 4.31909
4.6000 0.63182 4.36874
4.7000 0.63183 4.41771
4.8000 0.63183 4.46602
4.9000 0.63184 4.51371
5.0000 0.63184 4.56081
5.1000 0.63184 4.60733
5.2000 0.63184 4.65330
5.3000 0.63184 4.69874
5.4000 0.63184 4.74367
5.5000 0.63184 4.78811
5.6000 0.63185 4.83207
5.7000 0.63185 4.87558
5.8000 0.63185 4.91865
5.9000 0.63185 4.96129
6.0000 0.63185 5.00351
6.1000 0.63185 5.04534
6.2000 0.63185 5.08677
6.3000 0.63185 5.12783
6.4000 0.63185 5.16853
6.5000 0.63185 5.20886
6.6000 0.63185 5.24885
6.7000 0.63185 5.28851
6.8000 0.63185 5.32783
6.9000 0.63185 5.36684
7.0000 0.63185 5.40554
7.1000 0.63185 5.44393
7.2000 0.63185 5.48203
7.3000 0.63185 5.51984
7.4000 0.63185 5.55737
7.5000 0.63185 5.59462
7.6000 0.63185 5.63161
7.7000 0.63185 5.66833
7.8000 0.63185 5.70479
7.9000 0.63185 5.74100
8.0000 0.63185 5.77697
8.1000 0.63185 5.81269
8.2000 0.63185 5.84818
8.3000 0.63185 5.88344
8.4000 0.63185 5.91847
8.5000 0.63185 5.95328
8.6000 0.63185 5.98788
8.7000 0.63185 6.02225
8.8000 0.63185 6.05642
8.9000 0.63185 6.09039
9.0000 0.63185 6.12415
9.1000 0.63185 6.15771
9.2000 0.63185 6.19108
9.3000 0.63185 6.22426
9.4000 0.63185 6.25725
9.5000 0.63185 6.29006
9.6000 0.63185 6.32268
9.7000 0.63185 6.35513
9.8000 0.63185 6.38740
9.9000 0.63185 6.41950
10.0000 0.63185 6.45143
10.1000 0.63185 6.48319
10.2000 0.63185 6.51479
10.3000 0.63185 6.54623
10.4000 0.63185 6.57751
10.5000 0.63185 6.60863
10.6000 0.63185 6.63959
10.7000 0.63185 6.67041
10.8000 0.63185 6.70108
10.9000 0.63185 6.73159
11.0000 0.63185 6.76197
11.1000 0.63185 6.79220
11.2000 0.63185 6.82229
11.3000 0.63185 6.85224
11.4000 0.63185 6.88205
11.5000 0.63185 6.91173
11.6000 0.63185 6.94127
11.7000 0.63185 6.97068
11.8000 0.63185 6.99997
11.9000 0.63185 7.02912
12.0000 0.63185 7.05815
}\untw
\begin{tikzpicture}
\begin{axis}[width=9cm,height=6.5cm,xmin=0,xmax=12,ymin=0,ymax=7.5,
  axis lines=left,xlabel={$t$},ylabel={$p(t),\,q(t)$},ytick={2,4,6},
  title={},
  legend pos=north west,legend cell align=left,legend style={draw=none,font=\footnotesize,fill=none},
  tick label style={font=\footnotesize}]
\addplot[green!45!black,very thick] table[x=t,y=q]{\degone};
\addlegendentry{$\mathcal O(-1)$}
\addplot[fib,very thick] table[x=t,y=q]{\degtwo};
\addlegendentry{$\mathcal O(-2)$ (K\"ahler)}
\addplot[base,very thick,mark=o,mark size=1.3pt,mark repeat=10,mark phase=6,mark options={solid,fill=white}] table[x=t,y=q]{\untw};
\addlegendentry{$\mathbb S^2\times \mathbb R^2$ (untwisted)}
\addplot[green!45!black,very thick,dashed,forget plot] table[x=t,y=p]{\degone};
\addplot[fib,very thick,dashed,forget plot] table[x=t,y=p]{\degtwo};
\addplot[base,very thick,dashed,mark=o,mark size=1.3pt,mark repeat=10,mark phase=6,mark options={solid,fill=white},forget plot] table[x=t,y=p]{\untw};
\end{axis}
\end{tikzpicture}
\setlength{\abovecaptionskip}{1pt}
 \caption{\small Bolts with $\mu=-2$:
    $q(t)$ is solid and $p(t)$ is dashed.}
  \label{fig:boltprofiles}
\end{figure}

To illustrate these normalizations, Figure \ref{fig:boltprofiles} plots
numerical bolt profiles with $\mu=-2$ in twisted degrees one and two and
in the untwisted case.

\begin{remark}
    Distinct normalized parameters give non-homothetic metrics, and different rows or different bolt degrees are not identified. The nut $\delta=0$ is Bryant's soliton. The higher-degree endpoint $\mu_m$ is non-collapsed; every bolt with a smaller parameter is collapsed. Here non-collapse has its usual curvature-scale meaning: there is a positive lower bound for $\operatorname{Vol}B(x,r)/r^4$ whenever $\sup_{B(x,r)}|\operatorname{Rm}|\leq r^{-2}$. 
\end{remark}

\begin{remark}
For the natural complex structure, the K\"{a}hler member among the nuts is $\delta=-1/2$, which is H.-D. Cao's soliton \cite{Cao}. For degree two, every bolt with $\mu<0$ is the complete K\"{a}hler trajectory constructed by A. Dancer and M. Wang \cite{DW} and H.-D.~Cao~\cite{Cao}. 
Our analysis of the nut case recovers the existence part of \cite[Theorem 8.1]{Appleton}, cases (3) and (4); the exclusion of $\delta>0$ is not treated there. Furthermore, the higher-degree uniqueness in Theorem~\ref{threshold} (Section \ref{thresholdsec}) resolves the real 4D biaxial case of \cite[Conjecture 1]{Appleton}. 
\end{remark}

Our approach is to investigate the Weyl tensor, $W_g$, of a conformally K\"ahler metric under the leafwise transnormal hypothesis,
in the spirit of \cite{Derdzinski, CH2016}. Indeed, Theorem \ref{alignment} applies to any sign of $\lambda$ and will be of independent interest. Then, in the steady case, the problem reduces to the case of cohomogeneity one symmetry with specific isometry groups (Proposition \ref{globalsym}). Here the existence is already provided by \cite{Appleton, Ivey, Stolarski, DW, Wink}. Our contribution is to tie up some loose ends, identify appropriate parameters, and prove uniqueness. Roughly speaking, the argument is purely ODE-based, combining ideas from ~\cite{buttsworth21su2, Appleton, Donovan, TY, DW} with some new ingredients such as the integrating-factor identities (\ref{pweight})--(\ref{qweight}), the linearized system (\ref{firstorder}), and the barrier in Lemma \ref{barrier}.

The first step is to show that the conformal factor $u$ is
aligned with $f$, that is, $du\wedge df\equiv0$ thanks to an
identity for the Cotton tensor. We obtain the following
dichotomy for GRS of every sign, compact or non-compact, complete or incomplete.

\begin{theorem}\label{alignment}
Let $(M^4,\,g,\,f)$ be a connected GRS with non-constant leafwise transnormal potential $f$ whose metric is conformally K\"ahler. Then
\begin{equation}\label{dichotomy}
 W_g\equiv0\quad\text{or}\quad du\wedge df\equiv0.
\end{equation}
Every point with $du\wedge df\neq0$ has a neighborhood on which $g$ is an interval warped over a three-dimensional space form. 
\end{theorem}

Theorem~\ref{alignment} has two immediate consequences. The first is
that, away from the locally conformally flat case, alignment holds
unconditionally.

\begin{corollary}\label{alignmentnonlcf}
Under the hypotheses of Theorem~\ref{alignment}, if $g$ is not locally conformally flat, then $du\wedge df\equiv0$. In particular, near every regular point of $f$ there exists a smooth function $U$ of one variable such that $u=U(f).$ 
\end{corollary}

On a closed manifold, the locally conformally flat alternative never
occurs, and alignment holds automatically.

\begin{corollary}\label{alignmentclosed}
Under the hypotheses of Theorem~\ref{alignment}, if, in addition, $M$ is closed, then $\lambda>0$, $g$ is not locally conformally flat, and $du\wedge df\equiv0$ on $M$.
\end{corollary}

\begin{remark}\label{maschlerremark}
The alternative $W_g\equiv0$ cannot be removed from the local statement; see Remark~\ref{lcfremark}. Alignment is the functional-dependence hypothesis under which K\"{a}hler metrics conformal to GRS are studied in \cite[Section 6]{Maschler}. Proposition~\ref{normalform} adds the local normal form for every $\lambda$ under the leafwise transnormal hypothesis.
\end{remark}

The rest of the article is organized as follows. In
Section~\ref{preliminary}, we collect the useful Ricci soliton identities and scalar
relations used throughout the paper. Section~\ref{localsym} establishes
the alignment dichotomy of Theorem~\ref{alignment}, treats the Bryant
soliton in Proposition~\ref{bryant}, and derives the aligned local normal
form and its symmetry in Proposition~\ref{normalform}. In
Section~\ref{globalsec}, we globalize these symmetries and identify the
possible singular orbits and underlying manifolds (Proposition \ref{globalsym}).
Section~\ref{trajectories} analyzes the complete nut and bolt trajectories
and their asymptotic behavior. Section~\ref{thresholdsec} proves
Theorem~\ref{threshold}, establishing uniqueness of the higher-degree
shooting threshold and the associated non-collapsed bolt. Finally,
Section~\ref{converse} verifies the converse construction, establishes
irreducibility, distinguishes the normalized parameters, and completes
the proof of Theorem~\ref{classification}.

\textbf{Declaration of AI Usage:} This work was supported in part by a grant of access to OpenAI models through the ChatGPT for Academic Researchers program. All proofs are completely written by the authors, who retain exclusive responsibility for the content of this work.

\section{Preliminaries}\label{preliminary}

Throughout this section and the next, $(M^4,g,f)$ is a connected 4D GRS with non-constant leafwise transnormal potential. As in~\eqref{ck}, $k=e^{2u}g$ is K\"ahler with respect to $J$. The steady case is $\lambda=0$, and all one-variable data are local on regular foliation charts.
Our curvature convention is $R_{ijk\ell}=g(R(e_i,e_j)e_\ell,e_k)$ so that $R_{ijij}>0$ on a round sphere, and $R$ denotes the scalar curvature. We now collect several well-known identities.

\begin{lemma}[\cite{hamilton95formation}]
\label{lemmaA}
We have
\begin{equation}
\label{eq1a}
    R + \Delta f = 4\lambda,
\end{equation}
\begin{equation}
\label{eq2a}
    \nabla_i R=2 R_{ij}\nabla_j f,
\end{equation}
\begin{equation}
\label{eq3a}
      R+|\nabla f|^2 =2\lambda f + C,\,\,\,\,\hbox{for some constant $C$},
\end{equation}
\begin{equation}
\label{eq4a}
    \Delta_fR=2\lambda R-2|\Ric|^2,
\end{equation}
\begin{equation}
\label{eq5a}
    \nabla_l R_{ijkl}=\nabla_j R_{ik}-\nabla_{i}R_{jk}=R_{ijkl}\nabla_l f.
\end{equation} Here, $\Delta_f=\Delta-\langle\nabla f,\nabla\cdot\rangle$ denotes the drifted Laplacian operator. 
\end{lemma}

Proceeding, we record the Ricci soliton identities that we need on a regular foliation chart. By the leafwise transnormal hypothesis, on such a chart we write
\[
 |\nabla f|^2=\mathcal B(f)>0,\qquad
 \mathcal F=\sqrt{\mathcal B},\quad \nu=\nabla f/\mathcal F.
\]
Primes on functions of $f$ in this section and the next mean differentiation with respect to $f$. For $N_{\rm Ric}=|\Ric_g|^2$, by Lemma \ref{lemmaA}, one obtains that 
\begin{equation}\label{scalardata}
 \begin{gathered}
\Delta f=4\lambda-R,\\ R=C+2\lambda f-\mathcal B(f),\\ 
 r_0:=\Ric(\nu,\nu)=\tfrac12R'=\lambda-\tfrac12\mathcal B',\\
 \Ric(\nu,X)=0\quad \hbox{for}\quad X\perp\nu,\quad \nabla_\nu\nu=0,\\
 N_{\rm Ric}=\lambda R-\tfrac12\bigl[\mathcal B R''+(\Delta f-\mathcal B)R'\bigr].
 \end{gathered}
\end{equation}
Thus $R$, $\Delta f$, $r_0$, and $N_{\rm Ric}$ are functions of $f$. In particular, for every $\lambda$, the leafwise transnormal hypothesis implies $\Delta f$ is also constant on each connected component of a regular level set; compare \cite[Proposition 2.11]{Tran}.

When $\lambda=0$, we have \[\Delta f=\mathcal B-C,\quad \,r_0=-\tfrac12\mathcal B'\quad\hbox{and}\quad N_{\rm Ric}=\tfrac12(\mathcal B\mathcal B''-C\mathcal B').\] 
For the complex structure $J$ we set $$\omega(X,Y)=g(JX,Y).$$

We define the \textit{Cotton tensor} by
\begin{equation}\label{cotton}
 \begin{aligned}
 C_{ijk}&=\nabla_i\Ric_{jk}-\nabla_j\Ric_{ik}
           -\tfrac16(R_i g_{jk}-R_j g_{ik})\\&=-R_{ijk\ell}f_\ell-\tfrac16(R_i g_{jk}-R_j g_{ik}),\\     C_{ijk}&=-2\nabla^\ell W_{ijk\ell},
 \end{aligned} 
\end{equation} where we have used Eq. (\ref{eq5a}). Here $R_i$ and $f_\ell$ denote the components of $dR$ and $df$, respectively, and $W$ is the Weyl tensor of $g$. 

Because $dR$ is proportional to $df$, these identities imply
\begin{equation}\label{radialcotton}
 C(X,Y,\nu)=0\quad\text{for all }X,Y\in \mathfrak{X}(M).
\end{equation}

\section{Alignment and local symmetry}\label{localsym}

In this section, we restrict ourselves to dimension four and prove Theorem~\ref{alignment}, Corollaries~\ref{alignmentnonlcf} and \ref{alignmentclosed}, and the local normal form which is valid for any sign of $\lambda$. We also treat the conformally flat steady branch. 

\subsection{The K\"{a}hler Weyl calculation and tangential Ricci spectrum}
The K\"{a}hler Weyl formula \cite[Proposition 2]{Derdzinski} and conformal covariance give
\begin{equation}\label{kahlerweyl}
 W^+_{ijk\ell}=\mathfrak p K_{ijk\ell},\qquad
 K_{ijk\ell}=\tfrac32\omega_{ij}\omega_{k\ell}
 -\tfrac12(g_{ik}g_{j\ell}-g_{i\ell}g_{jk}+\varepsilon_{ijk\ell}),
\end{equation}
where $W^+$ is the self-dual part of the Weyl tensor of $g$, $\mathfrak{p}=\frac1{12}e^{2u}R_{e^{2u}g}$, and $\varepsilon$ is the volume form of $g$. On self-dual two-forms, the eigenvalues of $W^+$ are $2\mathfrak p,-\mathfrak p,-\mathfrak p$.

The next lemma computes the divergence of $W^+$. 

\begin{lemma}\label{weyldiv}
For the one-form $\beta_W=d\mathfrak p-3\mathfrak p\,du$, we have
\begin{equation}\label{weyldivergence}
 \nabla^\ell W^+_{ijk\ell}=K_{ijk\ell}(\beta_W)_\ell.
\end{equation}
\end{lemma}

\begin{proof}
This follows from the proof of \cite[Proposition 3]{Derdzinski}. The key is the conformal change formula for the divergence of the Weyl tensor \cite[Equation (19)]{Derdzinski}.

\end{proof}

For later use, fix an adapted oriented orthonormal frame
\[
e_1=Je_0,\qquad e_3=Je_2,
\qquad
\omega_{01}=\omega_{23}=1.
\]
Then, for every one-form \(z\),
\begin{equation}
\label{rem4}
 \bigl(K_{0a0\ell}z_\ell\bigr)_{a=1}^3 = \left(z_1,-\frac12z_2,-\frac12z_3\right),~  \bigl(K_{0ba0}\bigr)_{b,a=1}^3
 = -\operatorname{diag} \left(1,-\frac12,-\frac12\right).
\end{equation}
In particular, the map
\[
(z_1,z_2,z_3)
\longmapsto
\bigl(K_{0a0\ell}z_\ell\bigr)_{a=1}^3
\]
is invertible. These formulas will be used below in the Cotton-tensor
calculation and also fix the relevant signs.

We now work on the regular set and choose an adapted orthonormal frame
\[
e_0=\nu,\qquad e_1=\xi=J\nu,\qquad e_3=Je_2,
\]
with $\mathcal H=\operatorname{span}\{\nu,\xi\}^{\perp}.$ Recall that
\[
C^+_{ijk}=-2\nabla^\ell W^+_{ijk\ell}.
\]
Hence, by Lemma~\ref{weyldiv}, we have
\begin{equation}\label{eq:Cplus-beta}
 C^+_{ijk}=-2K_{ijk\ell}(\beta_W)_\ell .
\end{equation}

On the other hand, \eqref{radialcotton} gives $C(X,Y,\nu)=0,$ for all $X,Y\in \mathfrak{X}(M)$. Therefore, its self-dual projection in the first two slots also satisfies $C^+_{0a0}=0,$ $ a=1,2,3.$ Using \eqref{eq:Cplus-beta} yields
\[
K_{0a0\ell}(\beta_W)_\ell=0,
\qquad a=1,2,3.
\]
By the contraction formula in (\ref{rem4}), the
tangential part of \(\beta_W\) vanishes:
\begin{equation}\label{eq:tangential-beta}
 (\beta_W)|_{T\Sigma}=0 \iff d_\Sigma\mathfrak p= 3\mathfrak p\,d_\Sigma u.
\end{equation}
Next, since
\[
d(e^{-3u}\mathfrak p)
=
e^{-3u}(d\mathfrak p-3\mathfrak p\,du)
=
e^{-3u}\beta_W,
\]
equation \eqref{eq:tangential-beta} says that \(e^{-3u}\mathfrak p\)
is constant on each connected regular component. After restricting to a foliation
chart, there is a smooth function \(\Pi\) such that
\begin{equation}\label{eq:p-Pi}
 e^{-3u}\mathfrak p=\Pi(f).
\end{equation}
Differentiating \eqref{eq:p-Pi} gives
\begin{equation}\label{eq:beta-Pi}
 \beta_W=e^{3u}\Pi'(f)\,df.
\end{equation}
In particular, if $\beta_0:=\beta_W(\nu),$ then
\begin{equation}\label{eq:beta0}
 \beta_0=e^{3u}\Pi'(f)\mathcal F.
\end{equation}

Notice that no division by \(\mathfrak p\) has been used; consequently,
these identities remain valid through zeros of \(\mathfrak p\). The next result extracts the tangential Ricci spectrum from the remaining
Cotton components. We define the \emph{Ricci gap}
\[
\tau:=r_1-r_H.
\]

\begin{proposition}\label{spectrum}
On the regular set, we have, for $\Ric$ denoting the endomorphism,
\[
\Ric \nu=r_0\nu,\qquad
\Ric \xi=r_1\xi,\qquad
\Ric X=r_HX
\quad\text{for every }X\in\mathcal H.
\]
Thus \(\xi\) is a Ricci eigenvector, \(\mathcal H\) is a repeated Ricci
eigenspace, and all radial--tangential mixed components vanish. Locally on each regular foliation chart,  for a smooth function \(\Theta\),
\begin{equation}\label{eq:tau-level}
 \tau^2=\Theta(f) \quad\text{ and }\quad d_\Sigma\tau=3\tau\,d_\Sigma u.
\end{equation}
Consequently,
\begin{equation}\label{eq:compatibility-tau}
 \tau^2(du\wedge df)=0.
\end{equation}
All these identities remain valid through zeros of \(\mathfrak p\) and
\(\tau\). 
\end{proposition}

\begin{proof}
The radial--tangential mixed terms have already vanished by
\eqref{scalardata}. It remains to determine the tangential Ricci tensor.
We compute the self-dual Cotton component \(C^+_{0ba}\) in two different
ways.

Let
\[
S=(\Ric_{ab})_{a,b=1}^3,
\qquad
\mathsf D
=
\operatorname{diag}
\left(1,-\frac12,-\frac12\right).
\]
We use the four-dimensional curvature decomposition
\begin{equation}\label{eq:curv-decomp-spectrum}
 \begin{aligned}
 R_{ijk\ell}
 ={}&
 W_{ijk\ell}
 +\frac12\bigl(
 g_{ik}\Ric_{j\ell}
 -g_{i\ell}\Ric_{jk}
 -g_{jk}\Ric_{i\ell}
 +g_{j\ell}\Ric_{ik}
 \bigr)\\
 &-\frac R6
 \bigl(g_{ik}g_{j\ell}-g_{i\ell}g_{jk}\bigr).
 \end{aligned}
\end{equation}
Taking into account that
\[
\nabla f=\mathcal F e_0,
\qquad
\Ric(e_0,e_a)=0,
\qquad
R_0=2r_0\mathcal F,
\qquad
R_a=0,
\]
the Cotton identity \eqref{cotton} gives, for \(a,b\in\{1,2,3\}\),
\begin{equation}\label{eq:C0ba-unprojected}
 C_{0ba}
 =
 -\mathcal F W_{0ba0}
 +\frac{\mathcal F}{2}
 \left(
 S_{ab}-\frac{R-r_0}{3}\delta_{ab}
 \right).
\end{equation}
Indeed, substituting \eqref{eq:curv-decomp-spectrum} into
\[
C_{0ba}
=
-\mathcal F R_{0ba0}
-\frac16R_0\delta_{ab}
\]
gives exactly \eqref{eq:C0ba-unprojected}.

For purely tangential \(c,d,a\), the same Cotton formula gives
\[
C_{cda}=-\mathcal F R_{cda0}.
\]
Next, because all radial--tangential Ricci components vanish, the Ricci and scalar
parts of \eqref{eq:curv-decomp-spectrum} make no contribution to
\(R_{cda0}\). Hence
\[
C_{cda}=-\mathcal F W_{cda0}.
\]
Projecting the first two indices onto the self-dual part, we get
\[
T^+_{ijk}
=
\frac12
\left(
T_{ijk}
+\frac12\varepsilon_{ijcd}T_{cdk}
\right).
\]
Consequently,
\[
\begin{aligned}
 C^+_{0ba}
 &=
 -\frac{\mathcal F}{2}
 \left(
 W_{0ba0}
 +\frac12\varepsilon_{0bcd}W_{cda0}
 \right)
 +\frac{\mathcal F}{4}
 \left(
 S_{ab}-\frac{R-r_0}{3}\delta_{ab}
 \right)\\
 &=
 -\mathcal F W^+_{0ba0}
 +\frac{\mathcal F}{4}
 \left(
 S_{ab}-\frac{R-r_0}{3}\delta_{ab}
 \right).
\end{aligned}
\]
Since \(W^+=\mathfrak p K\), we infer from~(\ref{rem4}) that $K_{0ba0}=-\mathsf D_{ab}$. Thus, this becomes
\begin{equation}\label{eq:Cplus-spectrum}
 C^+_{0ba}
 =
 \mathcal F
 \left[
 \mathfrak p\,\mathsf D_{ab}
 +\frac14
 \left(
 S_{ab}-\frac{R-r_0}{3}\delta_{ab}
 \right)
 \right].
\end{equation}

We now compute the same component using Lemma~\ref{weyldiv}.
Since \(\beta_W\) has no tangential component by
\eqref{eq:beta-Pi},
\[
C^+_{0ba}
=
-2K_{0ba0}\beta_0
=
2\beta_0\,\mathsf D_{ab}.
\]
Comparing this with \eqref{eq:Cplus-spectrum} gives
\begin{equation}\label{riccigap}
 \begin{gathered}
 S-\frac{R-r_0}{3}\operatorname{Id}
 =
 \left(
 \frac{8\beta_0}{\mathcal F}
 -4\mathfrak p
 \right)\mathsf D,\\[2mm]
 r_1=\frac{R-r_0+2\tau}{3},
 \qquad
 r_H=\frac{R-r_0-\tau}{3},\\[1mm]
 \tau
 =
 \frac{12\beta_0}{\mathcal F}
 -6\mathfrak p
 =
 6e^{3u}\bigl(2\Pi'(f)-\Pi(f)\bigr).
 \end{gathered}
\end{equation}
Thus \(\xi\) is an eigenvector of \(\Ric\) and
\(\mathcal H\) is a repeated eigenspace. 
Next,
\[
N_{\rm Ric}
=
r_0^2+r_1^2+2r_H^2
=
r_0^2
+\frac{(R-r_0)^2}{3}
+\frac{2\tau^2}{3}.
\]
Therefore
\begin{equation}\label{eq:tau-square}
 \tau^2
 =
 \frac32\bigl(N_{\rm Ric}-r_0^2\bigr)
 -\frac12(R-r_0)^2.
\end{equation}
By \eqref{scalardata}, every quantity on the right-hand side is a function
of \(f\). Hence \(\tau^2=\Theta(f)\) locally on each foliation chart.

Finally, differentiating the last expression for \(\tau\) in
\eqref{riccigap} tangentially yields
\[
d_\Sigma\tau
=
3\tau\,d_\Sigma u,
\]
because \(\Pi\) and \(\Pi'\) depend only on \(f\), which is constant along
each regular level. This jointly with the fact that
\(\tau^2=\Theta(f)\) gives
\[
0=d_\Sigma(\tau^2)
=
2\tau\,d_\Sigma\tau
=
6\tau^2d_\Sigma u.
\]
Consequently,
\begin{equation}\label{eq:tau-dSigma-u}
 \tau^2d_\Sigma u=0.
\end{equation}
Since \(df\) is normal to the regular levels, only the tangential component
of \(du\) contributes to \(du\wedge df\). Equation
\eqref{eq:tau-dSigma-u} is therefore equivalent to $\tau^2(du\wedge df)=0.$
\end{proof}


\subsection{The exceptional umbilic branch and global dichotomy}
We are ready to give the proof of Theorem~\ref{alignment}.

\begin{proof}[{\bf Proof of Theorem~\ref{alignment}}]
To begin with, we set
\[
\mathcal O:=\{x\in M:(du\wedge df)_x\neq0\}.
\]
Since $du\wedge df=0$ whenever $df=0$, one sees that $
\mathcal O\subset\mathcal R:=\{df\neq0\}.$ By Proposition~\ref{spectrum}, we have $\tau^2(du\wedge df)=0$ on $\mathcal R.$ Hence
\begin{equation}\label{eq:tau-zero-O}
 \tau=0\qquad\text{on }\mathcal O.
\end{equation}

We now divide the proof into two steps. 

\medskip
\noindent
\textbf{Step 1: the umbilic branch.}
Let $\mathcal U\subset\mathcal R$ be an open set on which $\tau=0$.
Then Proposition~\ref{spectrum} gives
\[
\Ric|_{T\Sigma}=r_H(f)\operatorname{Id},
\qquad
r_H=\frac{R-r_0}{3},
\]
where $\Sigma$ denotes a regular level of $f$. Moreover, taking into account that $R$ and $r_0$
are functions of $f$, so is $r_H$.

Since $\mathcal F=|\nabla f|$ and $\nu=\nabla f/\mathcal F,$ by leafwise transnormality, one has $\mathcal F=\mathcal F(f)$, and therefore
$\nu$ is geodesic. For tangent vectors $X,Y\in T\Sigma$,
\[
\Hess f(X,Y)
=
\mathcal F\,g(\mathsf L X,Y),
\]
where $\mathsf L X=\nabla_X\nu$ is the shape operator of the level set. By (\ref{grs}), we get
\[
\mathsf L
=
\frac{1}{\mathcal F}
\left(\lambda\operatorname{Id}
      -\Ric|_{T\Sigma}\right)
=
\frac{\lambda-r_H(f)}{\mathcal F(f)}
\operatorname{Id}.
\]
Thus every regular level in $\mathcal U$ is totally umbilic, with principal
curvature depending only on $f$. Choose a distance coordinate $s$ along the integral curves of $\nu$, so
that
\[
\partial_s=\nu,\qquad
g=ds^2+g_s,\qquad
f=f(s).
\]
Since $\frac12\partial_s g_s(X,Y)=g(\mathsf L X,Y),$ we obtain
\[
\partial_s g_s
=
2\frac{\lambda-r_H(f(s))}{\mathcal F(f(s))}\,g_s.
\]
The coefficient depends only on $s$. Integrating this scalar evolution yields
\begin{equation}\label{eq:umbilic-warped}
 g=ds^2+\psi(s)^2h_3,
\end{equation}
where $h_3$ is independent of $s$. In particular, for a warped product of the form \eqref{eq:umbilic-warped}, the
tangential Ricci tensor is
\[
\Ric_g|_{T\Sigma}
=
\Ric_{h_3}
-
\bigl(\psi\psi''+2(\psi')^2\bigr)h_3.
\]
Since
\[
\Ric_g|_{T\Sigma}
=
r_H(f(s))\,g_s
=
r_H(f(s))\psi^2h_3,
\]
we infer
\begin{equation}\label{eq:h3-Einstein}
 \Ric_{h_3}
 =
 \left[
 \psi\psi''+2(\psi')^2
 +r_H(f(s))\psi^2
 \right]h_3.
\end{equation}
Fixing any value of $s$ shows that $h_3$ is Einstein. This therefore implies that $h_3$ is a
three-dimensional space form.

Consequently, \(g\) is locally conformally flat on \(\mathcal U\).
Indeed, after introducing a coordinate \(\rho\) by $d\rho=\frac{ds}{\psi(s)},$
we may write $g=\psi(s)^2\bigl(d\rho^2+h_3\bigr),$
and the product of an interval with a three-dimensional space form is
conformally flat. Equivalently, the warped-product curvature formulas give
directly
\[
W_g=0
\qquad\text{on }\mathcal U.
\]

\medskip
\noindent
\textbf{Step 2: the global dichotomy.}
Suppose that $\mathcal O\neq\varnothing$. By
\eqref{eq:tau-zero-O}, the Ricci gap vanishes throughout the open set
$\mathcal O$. Step~1 therefore gives $W_g=0\text{ on }\mathcal O.$ Since every Ricci soliton metric is real analytic (see \cite[Corollary~1.3]{Kotschwar}) and $M$ is connected, by analytic continuation, one deduces that $W_g\equiv0
\text{ on }M.$

Finally, if $x\in M$ satisfies $(du\wedge df)_x\neq0$, then, it holds throughout an open neighborhood. Hence $\tau=0$ there
by Proposition~\ref{spectrum}, and Step~1 gives the asserted local warped
product over a three-dimensional space form. 
\end{proof}


\begin{proof}[{\bf Proof of Corollary~\ref{alignmentnonlcf}}]
Since $W_g\not\equiv0$, the second alternative in (\ref{dichotomy}) holds. On a regular foliation chart, $du\wedge df=0$ makes $u$ constant on the connected local levels, and a transverse section shows that $U$ is smooth.
\end{proof}

\begin{proof}[{\bf Proof of Corollary~\ref{alignmentclosed}}]
We necessarily have $\lambda>0$ and closed locally conformally flat shrinking GRS of dimension at least four are trivial \cite[Proposition 3.8]{ELM}. Thereby, since $f$ is non-constant, we have $W_g\not\equiv0$. Corollary~\ref{alignmentnonlcf} then applies.
\end{proof}

\begin{remark}\label{closeddirect}
The conformally K\"{a}hler structure also excludes the conformally flat alternative on a closed shrinking Ricci soliton directly, even when $f$ is constant. Indeed, if $W_g\equiv0$, then $W^+=0$, so $R_{\k}=0$ by (\ref{kahlerweyl}). Thus $v=e^u$ satisfies
\begin{equation}
\label{lk127}
 -6\Delta_gv+Rv=R_{\k}\,v^3=0.
\end{equation}
It is known from \cite{Chen} that every shrinking Ricci soliton has $R\geq0.$ Moreover, $R$ is strictly positive unless $(M^4,\,g,\, f)$ is the Gaussian shrinking soliton. Therefore, multiplying the conformal scalar-curvature equation (\ref{lk127}) by $v$ and integrating, we obtain
\[
0
=
\int_M
\left(
6|\nabla v|^2+Rv^2
\right)\,d\mathrm{vol}_g.
\] The integrand is nonnegative, while $v>0$ and $R\not\equiv0$, so the right-hand side is strictly positive, a contradiction.
\end{remark}

\begin{remark}\label{lcfremark}
The locally conformally flat case is well understood in the shrinking and steady settings, while in the expanding case the complete classification remains open \cite{CMMR, CC, CM, Zhang}. 
This explains the non-locally-conformally-flat hypothesis in Corollary~\ref{alignmentnonlcf}. Additionally, the conformally flat alternative in Theorem~\ref{alignment} cannot be removed from the local statement. Indeed, for $a\in\R^4$, let $g=\delta$ on $\R^4\setminus\{a\}$, with $f=\lambda|x|^2/2$ if $\lambda\neq0$ and $f=x_1$ if $\lambda=0$; for $\lambda\neq0$, take $a\neq0$. These are leafwise transnormal solitons of the given sign. The metric $k=|x-a|^{-4}\delta$ is the pull-back of the flat K\"{a}hler metric under the inversion centered at $a$, so it is K\"{a}hler for the pulled-back complex structure. But $u=-2\log|x-a|$ has $du\wedge df\not\equiv0$. 
\end{remark}

\subsection{The conformally flat steady branch}\label{lcfsteady}
In this subsection we assume completeness and $\lambda=0$. H.-D. Cao and Q. Chen~\cite[Theorem~1.2]{CC} proved that every complete, non-flat, locally conformally flat steady GRS of dimension at least three is, up to scaling, the Bryant soliton. In our setting, the flat case is excluded. It remains, however, to verify the alignment condition for an arbitrary globally defined K\"ahler representative conformal to the Bryant soliton; this does not follow directly from the rotational symmetry of the soliton metric $g$.

\begin{proposition}\label{bryant}
Let $(\mathbb R^4,g,f)$ be the 4D Bryant soliton. Suppose that $g$ is conformal to a K\"{a}hler metric $e^{2u}g$ (\ref{ck}). Then $u$ is radial. In particular, $du\wedge df=0$.
\end{proposition}

\begin{proof}
Write the Bryant metric in its rotationally symmetric form
\[
g=dt^2+\psi(t)^2g_{S^3},
\qquad t\in[0,\infty),
\]
where $t=0$ is the smooth pole. Thus $\psi(0)=0,$ $\psi'(0)=1,$ $
\psi(t)>0,$ for $t>0.$ The Bryant soliton~\cite{Bryant} has positive sectional curvature. In particular, its
radial sectional curvature satisfies $-\psi''/\psi >0,$ and hence $\psi''<0,$ for $t>0.$ Since $\psi'(0)=1$, concavity gives
\begin{equation*}
0<\psi(t)\le t.
\end{equation*}

Next, define a new radial coordinate $r$ on $(0,\infty)$ by
\begin{equation}\label{eq:bryant-r}
\frac{dr}{r}=\frac{dt}{\psi(t)}.
\end{equation} 
It is immediate to check that $r\to0$ as $t\to 0$ and $r(t)\longrightarrow\infty$ as $t\longrightarrow\infty$; thus $r$ gives a global radial coordinate on $\mathbb R^4$. 
From \eqref{eq:bryant-r},
\[
dt=\frac{\psi(t)}{r}\,dr,
\text{ and hence }
g
=
\frac{\psi(t(r))^2}{r^2}
\bigl(dr^2+r^2g_{S^3}\bigr)
=
a(r)^2\delta,
\]
where $a(r):=\frac{\psi(t(r))}{r}>0.$ The regularity of $\psi$ near the pole shows that $a$ extends smoothly and positively
across $r=0$. Thus the Bryant metric is globally conformal to the Euclidean
metric on $\mathbb R^4$.

Now set $v:=e^u a>0.$ This implies $e^{2u}g=v^2\delta.$ Next, since $g$ is conformally flat, so is $\k=e^{2u}g$. As it is
K\"ahler, its self-dual Weyl tensor is determined by its scalar curvature~\cite[Proposition 2]{Derdzinski}; hence, $W_{\k}^+=0$ implies $R_{\k}=0.$ The four-dimensional conformal
scalar-curvature formula therefore gives
\[
0=R_{\k}=-6v^{-3}\Delta_\delta v.
\] Thus $v$ is a positive harmonic function on all of $\mathbb R^4$, and
Liouville's theorem implies that $v\equiv c$ for some constant $c>0$.
Consequently, $u$ is radial.

Finally, the potential of the Bryant soliton is radial as well,
$f=f(t)=f(t(r))$. Hence both $du$ and $df$ are multiples of $dr$, and
therefore $du\wedge df=0,$ as asserted. 
\end{proof}

\begin{remark}
    We highlight that the proof uses only that $g=a(r)^2\delta$ globally on $\R^4$ and that $f$ is radial. Therefore, it applies verbatim to every rotationally symmetric soliton on $\R^4$ with radial potential and $\int^\infty_{t_0} \frac{dt}{\psi}=\infty$, $t_0 >0$, of any sign.
\end{remark}

\subsection{The aligned local normal form and its symmetry}
Here we show that alignment produces a K\"ahler Killing potential and, consequently, local cohomogeneity one symmetry. A foliation chart is called \emph{aligned} if $du\wedge df=0$ there.

\begin{proposition}\label{normalform}
Let $(M^4,g,f)$ be a GRS with non-constant leafwise transnormal potential, and suppose that $k=e^{2u}g$ is K\"ahler with respect to $J$. On an aligned regular foliation chart, there is a coordinate $s=s(f)$ with $s_f>0$ such that
\begin{equation}\label{calabi}
 e^{2u}g=\frac{ds^2}{V(s)}+V(s)\eta^2+\ell(s)h,\qquad
 \ell(s)=cs+d_0>0,\qquad d\eta=c\omega_h.
\end{equation}
Here $h$ is a surface metric of constant Gaussian curvature $\kappa$, $V>0$, and $u=v(s)$, $f=\phi(s)$. There is a four-dimensional local Lie algebra of Killing fields for both $g$ and $e^{2u}g$, preserving $f,u,J$, with evaluation of rank three tangent to the regular levels. This includes $c=0$ and every regular level where the Ricci gap $\tau$ vanishes.
\end{proposition}

\begin{proof}
The idea is to reparametrize $f$ so that its Hessian with
respect to the K\"ahler metric $\k=e^{2u}g$ is $J$-invariant. The resulting
Hamiltonian vector field is holomorphic and Killing, and its local quotient
produces the Calabi-type form \eqref{calabi}. 


Alignment gives $u=U(f)$ on the foliation chart. By \eqref{scalardata} and
\eqref{riccigap}, the Ricci eigenvalues $r_0,$ $r_1,$ $r_H$ are smooth functions of $f$. This remains true across zeros of the Ricci
gap $\tau=r_1-r_H$. Choose $s=s(f)$ by solving
\begin{equation}\label{reparam}
 \frac{s_{ff}}{s_f}
 =
 \frac{r_0-r_1}{\mathcal B}+2U_f,
 \qquad s_f>0.
\end{equation} Thus, since $\mathcal B=|\nabla f|^2>0$ on the regular chart, this is a smooth
scalar ODE and admits a local solution with $s_f>0$. We now claim that the Hessian of $s$ is $J$-invariant. Indeed, the conformal Hessian formula gives
\[
\Hess_{\k} f
=
\Hess_g f
-2U_f\,df\otimes df
+U_f\mathcal B\,g.
\] By (\ref{grs}), its diagonal entries in the adapted $g$-orthonormal frame $(\nu,\xi,e_2,e_3)$ are
\[
\lambda-r_0-U_f\mathcal B,\qquad
\lambda-r_1+U_f\mathcal B,\qquad
\lambda-r_H+U_f\mathcal B,\qquad
\lambda-r_H+U_f\mathcal B,
\]
and all mixed components vanish. Now, taking into account that 
\[
\Hess_{\k}s
=
s_f\Hess_{\k}f+s_{ff}\,df\otimes df,
\]
\eqref{reparam} gives
\[
\begin{aligned}
(\Hess_{\k}s)(\nu,\nu)-(\Hess_{\k}s)(\xi,\xi)
&=
s_f(r_1-r_0-2U_f\mathcal B)
+s_{ff}\mathcal B=0.
\end{aligned}
\]
Thus, $\Hess_{\k}s$ is $J$-invariant, as claimed. Its restriction to
\[
\mathcal H=\operatorname{span}\{\nabla^{\k}s,J\nabla^{\k}s\}^{\perp}
\]
has the form
\[
\Hess_{\k}s|_{\mathcal H}
=
\psi_H(s)\,\k|_{\mathcal H}
\]
for a smooth function $\psi_H$. Moreover,
\begin{equation}\label{eq:V-def}
 V(s):=|\nabla^{\k}s|_{\k}^2
 =
 e^{-2U}s_f^2\mathcal B>0.
\end{equation}
All these quantities depend only on $s$.

To proceed, denote $X:=\nabla^{\k}s$ and $K_0:=JX.$ Since $\k$ is K\"ahler and $\Hess_{\k}s$ is $J$-invariant, the Hessian endomorphism commutes with $J$. Consequently,
\[
\mathcal L_{K_0}\k=0
\qquad\hbox{and}\qquad
\mathcal L_{K_0}J=0,
\]
where $\mathcal L$ denotes the Lie derivative. Therefore $K_0$ is a holomorphic Killing field. The same commutation relation
also gives $\mathcal L_XJ=0.$ By denoting $N_s:=\frac{X}{V},$ one obtains that $N_s(s)=1.$ Since the Hessian endomorphism preserves
$\operatorname{span}\{X,K_0\}$ and has the same eigenvalue on $X$ and
$K_0$, one has $[X,K_0]=0.$ Furthermore \(K_0(V)=0\), because \(V=V(s)\). Consequently,
\begin{equation}\label{eq:NKcommute}
 [N_s,K_0]=0.
\end{equation}

Considering $\eta:=\frac{\k(K_0,\cdot)}{V}$, one deduces that
\[
\eta(K_0)=1,\qquad
\eta(N_s)=0,
\qquad
\mathcal H=\ker ds\cap\ker\eta.
\]
The normal flow preserves this splitting. Indeed, using
\eqref{eq:NKcommute} and the fact that the Hessian is scalar on
$\operatorname{span}\{X,K_0\}$, one obtains that
\begin{equation}\label{eq:eta-invariant}
 \mathcal L_{N_s}\eta=0.
\end{equation}
Moreover,
\begin{equation}\label{eq:horizontal-evolution}
 (\mathcal L_{N_s}\k)|_{\mathcal H}
 =
 \frac{2\psi_H(s)}{V(s)}\,\k|_{\mathcal H},
 \qquad
 (\mathcal L_{N_s}J)|_{\mathcal H}=0.
\end{equation}

Fix a reference level and use the flow of $N_s$ to identify nearby regular
levels. By \eqref{eq:eta-invariant}, the one-form $\eta$ is independent of
$s$ under this identification, while \eqref{eq:horizontal-evolution} shows that the horizontal metric evolves
only by a scalar factor depending on $s$. Hence there exist a fixed
surface metric $h$ and a positive function $\ell(s)$ such that \[\k|_{\mathcal H}=\ell(s)h.\] Since the horizontal complex structure is also preserved by the $N_s$-flow, if $\omega_h$ denotes the area form of $h$, then
\[
\Omega_{\k}|_{\mathcal H}
=
\ell(s)\omega_h.
\]

Next, taking into account that $\k(N_s,N_s)=\frac1V,$ $\k(K_0,K_0)=V,$ and $\k(N_s,K_0)=0,$ we infer
\begin{equation}\label{eq:pre-Calabi}
 \k=\frac{ds^2}{V}+V\eta^2+\ell(s)h,
 \qquad
 \Omega_{\k}=ds\wedge\eta+\ell(s)\omega_h.
\end{equation} In particular, \eqref{eq:eta-invariant} and Cartan's formula guarantee $\iota_{N_s}d\eta=0,$
while $\mathcal L_{K_0}\eta=0$ and $\eta(K_0)=1$ give $\iota_{K_0}d\eta=0.$ Thus \(d\eta\) is a horizontal two-form.

Now, taking the exterior derivative of the second identity in
\eqref{eq:pre-Calabi} and using $d\Omega_{\k}=0$ yields
\[
0
=
ds\wedge\bigl(\ell_s\omega_h-d\eta\bigr).
\]
Therefore $d\eta=\ell_s\,\omega_h.$ Both $\eta$ and $\omega_h$ are independent of $s$, so $\ell_s$ must be a
constant, say $c$. Consequently,
\[
\ell(s)=cs+d_0,
\qquad
d\eta=c\,\omega_h,
\]
which proves the stated Calabi form.

It remains to show that the Gaussian curvature of $h$ is constant. To this end, return
to $g=e^{-2U}\k$ and introduce the distance coordinate $t$ and the functions
\[
dt=\frac{e^{-U}}{\sqrt V}\,ds,
\qquad
p=e^{-U}\sqrt V,
\qquad
q=e^{-U}\sqrt\ell.
\]
Therefore 
\[
g=dt^2+p(t)^2\eta^2+q(t)^2h,
\qquad
d\eta=c\,\omega_h.
\]
Allow initially the Gaussian curvature $K_h$ of $h$ to depend on the
surface point. The horizontal component of the Ricci soliton equation is
\[
\lambda
=
\frac{K_h(x)}{q^2}
-\frac{c^2p^2}{2q^4}
-\frac{q''}{q}
-\frac{p'q'}{pq}
-\frac{(q')^2}{q^2}
+\frac{f'q'}q,
\]
Here primes denote $t$-derivatives. All terms other than $K_h(x)/q(t)^2$
depend only on $t$. Fixing $t$ therefore shows that
$K_h(x)$ is constant on the surface. Denote this constant by $\kappa$.

Finally, we shall establish the local symmetry algebra. An oriented surface of constant Gaussian curvature has a
three-dimensional local Lie algebra of Killing fields. Let $Y$ be one of
these fields. 

Since \(Y\) preserves the area form,
\[
d(\iota_Y\omega_h)
=
\mathcal L_Y\omega_h
=
0.
\]
Locally we may therefore choose a function $\varphi_Y$ satisfying $d\varphi_Y=-c\,\iota_Y\omega_h.$ Let \(Y^{\mathcal H}\) be the horizontal lift of \(Y\) with respect to
\(\eta\), and define $\widehat Y
=
Y^{\mathcal H}+\varphi_YK_0.$ Taking into account $\eta(\widehat Y)=\varphi_Y,$
Cartan's formula and \(d\eta=c\omega_h\) give
\[
\mathcal L_{\widehat Y}\eta
=
\iota_{\widehat Y}d\eta+d(\eta(\widehat Y))
=
c\,\iota_Y\omega_h+d\varphi_Y
=
0.
\]
The lifted field also preserves $s,h,\omega_h$, and consequently,
$\mathcal L_{\widehat Y}\k
=
\mathcal L_{\widehat Y}\Omega_{\k}
=
0.$ From this, it follows that $\widehat Y$ preserves \(J\). Since $u$ and $f$ depend
only on $s$, it also preserves $u,f$, and hence \[g=e^{-2u}\k.\]

Together with $K_0$, the lifts of the three local Killing fields of
$(N,h)$ form a four-dimensional local Lie algebra (compare \cite[Lemma 4.8]{tran23kahler}). To see closure, if
$Y_1,Y_2$ are base Killing fields, then
\[
[\widehat Y_1,\widehat Y_2]
-
\widehat{[Y_1,Y_2]}
\]
projects trivially to the quotient and is therefore of the form
\(aK_0\). Both terms preserve $\eta$, so
\[
0=\mathcal L_{aK_0}\eta=da,
\]
and \(a\) is constant. Moreover \(K_0\) commutes with all the lifted
fields. Thus the span is a four-dimensional Lie algebra.

At each point of a regular level, the evaluations of the base Killing
fields span the two horizontal directions, while $K_0$ spans the vertical
direction. Hence the evaluation of this Lie algebra has rank three and is
tangent to the regular levels. This proves the proposition.
\end{proof}



\section{Global symmetry and endpoints}\label{globalsec}

Throughout this section, we assume the hypotheses of Theorem~\ref{classification}. We globalize the symmetry of Proposition~\ref{normalform} and find the possible endpoints. We begin by recording the radial equations that will be used repeatedly. For a connection metric
\[
 g=dt^2+p^2\eta^2+q^2h,\qquad d\eta=c\omega_h,\quad K_h=\kappa,
\]
the steady soliton equation becomes
\begin{equation}\label{ode}
 \begin{aligned}
 f''&=\frac{p''}{p}+2\frac{q''}{q},\\
 p''&=\frac{c^2p^3}{2q^4}-2\frac{p'q'}q+f'p',\\
 q''&=\frac\kappa q-\frac{c^2p^2}{2q^3}
          -\frac{p'q'}p-\frac{(q')^2}q+f'q'.
 \end{aligned}
\end{equation}
From now on, primes denote differentiation with respect to $t$. Along any principal normal geodesic, set
\[
 F=f',\quad H=\operatorname{tr}\mathsf L,\quad
 D=\operatorname{tr}(\mathsf L^2)\geq H^2/3,\quad \Xi=H-F,
\]
where $\mathsf L$ is the shape operator of the principal orbits. For a connection metric, $H=p'/p+2q'/q$. The standard radial identities are
\begin{equation}\label{radialid}
 \begin{gathered}
 R+F^2=C,\quad \Xi'=-D,\quad F'+\Xi F=-C,\\
 H'=F'-D,\qquad R=-F'-HF,\qquad F''=DF-\Xi F'.
 \end{gathered}
\end{equation}

By \cite[Theorem 1.3]{ZhangPAMS}, \cite[Corollary 2.5]{Chen}, and 
\cite[Proposition 3.2]{PW09}, irreducibility gives $R>0$, hence $C>0$ and $|F|\leq\sqrt C$. We now turn to the global structure. The main result of this section is the following.

\begin{proposition}\label{globalsym}
There exists a compact connected four-dimensional group $G$ of isometries on $M$ preserving $f,u,$ and $J$ such that:
\begin{enumerate}[label=(\alph*)]
    \item Its orbit space is $[0,\infty)$, parametrized by the distance coordinate $t$, and its unique singular orbit $N=\Crit(f)$ is either a point or a round $\mathbb S^2$. 
    \item $M$ is the normal bundle of $N$ and its metric is the global ansatz (\ref{ansatz}).
    \item Every regular level of $f$ is a connected principal orbit. With $f|_N=0$, one has $f'<0$, $f''<0$, and $f(M)=(-\infty,0]$. Moreover, $|\nabla f|^2=B(f)$ and $u=U(f)$ for global functions smooth through the critical value. 
    \item Up to isomorphism, $G$ is $U(2)$ at a nut, $U(2)/\mathbb Z_m$ at a twisted bolt of degree $m$, and $SO(3)\times \mathbb S^1$ at an untwisted bolt. These groups act in the standard way on $\R^4=\C^2$, $\OO(-m)$, and $\mathbb S^2\times \mathbb R^2$, respectively. Here $\mathbb Z_m$ is the group of scalar matrices $\zeta I$ with $\zeta^m=1$.
\end{enumerate}
\end{proposition}

The proof occupies the rest of this section.

\subsection{Continuation of the symmetry algebra}
The idea is to extend the local Killing fields of Proposition~\ref{normalform} by analytic continuation, and then to check that the extensions preserve $f$, $u$, and $J$. GRS and their soliton vector fields are real analytic \cite{ivey96local, Kotschwar, DW}; in particular, $f$ is analytic. By Nomizu's continuation theorem~\cite{Nomizu} and the simple connectivity of $M$, each local Killing field of Proposition~\ref{normalform} extends uniquely to a global Killing field of $g$. By uniqueness, the bracket of two extensions is the extension of their bracket, so the extended fields span a four-dimensional Lie algebra $\mathfrak g$ of Killing fields.


We claim first that the extension preserves $f$. Let $X$ be one of the extended fields. Since $X$ is Killing, its local flow preserves the Levi-Civita connection and the Ricci tensor. Lie differentiating the steady Ricci soliton equation \eqref{grs} therefore gives
\[
\Hess_g(Xf)=0.
\] Hence, $\nabla(Xf)$ is parallel and $Xf$ vanishes on the original open set, so $Xf=0$ everywhere. 

By alignment (Theorem~\ref{alignment}, together with \cite[Theorem 1.2]{CC} and Proposition~\ref{bryant} in the conformally flat case), $Xu=0$ on the regular set. This set is dense; hence, $Xu=0$ on $M$. Thus $X$ is Killing for $\k=e^{2u}g$ as well.

The extension also preserves $J$, without any analyticity assumption on either $u$ or $J$. Indeed, $\nabla^{\k}J=0$, and the flow of a $\k$-Killing field preserves $\nabla^{\k}$. It follows that $\mathcal L_XJ$ is parallel. It vanishes on the original open set, hence vanishes globally. Completeness of $g$ makes each $X$ complete. 

Let $G$ be the closure in $\operatorname{Isom}(M,g)$ of the connected
subgroup with Lie algebra $\mathfrak g$. Then $G$ is a connected Lie
group acting properly on $M$. It preserves $f$, $u$ and $J$: the first
two properties pass to limits by continuity, and so does the third,
since the compact-open and $C^\infty$ topologies agree on
$\operatorname{Isom}(M,g)$.

We next determine its dimension. Let $p$ be a regular point of $f$. The isotropy group $G_p$ fixes
$\nu=\nabla f/|\nabla f|$, because it preserves $f$ and $g$, and it
fixes $J\nu$, because it preserves $J$. Hence it acts orthogonally and
$J$-linearly, that is, unitarily, on $\operatorname{span}\{\nu,J\nu\}^{\perp}$. An isometry of a connected
manifold is determined by its value and differential at one point, so
the isotropy representation is faithful and $\dim G_p\le1$. The orbit
$G\cdot p$ lies in the level set of $f$ through $p$, so
$\dim(G\cdot p)\le3$. Together with $\dim G\ge\dim\mathfrak g=4$, we obtain $\dim G=4$, $\dim G_p=1$ and $\dim(G\cdot p)=3$ at every
regular point. Thus $G$ acts properly with cohomogeneity one. In
particular, the subgroup generated by $\mathfrak g$ is open in the
connected group $G$, so it is already closed. 

\subsection{The orbit interval and the critical set}
\label{subsection4.3}
For a proper cohomogeneity-one action, the slice theorem of R.~Palais~\cite{Palais} implies that the orbit space is a line, a circle, a compact interval, or a half-line. Completeness of $g$ makes the quotient length metric complete. Thus, in the non-compact cases, the distance parameter ranges over all of $\R$ or $[0,\infty)$, respectively.
We first rule out the line case; compare \cite[Propositions 3.5 and 3.6]{TY}.

\begin{lemma}
\label{noline}
The orbit space is not homeomorphic to $\R$.
\end{lemma}

\begin{proof}
Suppose that $M/G\cong\R$. A unit-speed normal geodesic meets every orbit once. Along it, $F$, $H$, $D$
and $\Xi$ are smooth on $\R$ and satisfy (\ref{radialid}), with $C>0$, $D\ge H^2/3$ and $|F|<\sqrt C$, since $R=C-F^2>0$.
\begin{enumerate}
 \item[(i)] $|\Xi|\le\sqrt C$. Suppose $\Xi(t_1)<-\sqrt C$. For $t\ge t_1$ we have $\Xi\le\Xi(t_1)$, so
$H=\Xi+F<\Xi+\sqrt C<0$. Then $y:=-\Xi-\sqrt C$ satisfies $$y'=D\ge H^2/3\ge y^2/3$$ with $y(t_1)>0,$ so $y$ blows up
before $t_1+3/y(t_1),$ contradicting the smoothness of $\Xi$ on $\R.$ The case $\Xi(t_1)>\sqrt C$ is the same, since (\ref{radialid}) is invariant under
$(t,F,\Xi)\mapsto(-t,-F,-\Xi)$.

\item[(ii)] $\Xi$ is non-increasing. So the limits $\Xi_\pm:=\lim_{t\to\pm\infty}\Xi$ exist and $\Xi_-\ge\Xi_+$.

\item[(iii)] $\Xi_+=\sqrt C$. Otherwise fix $T$ with $\xi:=\Xi(T)<\sqrt C$. For $t\ge T$ we have
$-\sqrt C\le\Xi\le\xi$.
\begin{itemize}
  \item Where $F\ge0$, $F'=-C-\Xi F\le-C+\sqrt C\,F$. So $F$ decreases, at a rate bounded away from zero, until it
  becomes negative.
  \item Where $F<0$, $F'=-C+\Xi|F|\le-\sqrt C\,(\sqrt C-\max\{\xi,0\})<0$.
\end{itemize}
Hence $F\to-\infty$, which contradicts $|F|<\sqrt C$. By the same symmetry, $\Xi_-=-\sqrt C$.

\item[(iv)] Now $\Xi_-=-\sqrt C<\sqrt C=\Xi_+$, which contradicts (ii). 
\end{enumerate}
\end{proof}

We next analyze the boundary of the orbit space and the behavior of the potential. At an interior critical point of the descended potential, $F=0$, and~\eqref{radialid} gives $F'=-C<0$. Hence every such point is a strict local maximum.
At a boundary orbit, invariance under the slice isotropy implies that $df=0$ and that the restriction of $\Hess f$ to the normal space is $\mu_0\operatorname{Id}$. Since $f$ is constant along the orbit, the tangential part of its Hessian vanishes at the boundary. If the boundary orbit has codimension $d$, equation (\ref{eq1a}) gives
\begin{equation}\label{criticalhess}
 d\,\mu_0=\Delta f=-R.
\end{equation}
Since $F=0$ there, Hamilton's identity (\ref{eq3a}) gives $R=C$, and hence
\begin{equation}\label{criticalhess2}
 \mu_0=f''(0)=-\frac{R(0)}{d}
 =-\frac{C}{d}<0.
\end{equation}
Thus every boundary point is a strict one-sided local maximum of the potential. If the orbit space were a circle or a compact interval, the descended potential would attain a minimum, whereas every critical point, whether interior or boundary, is a strict local maximum. Consequently, the orbit space cannot be compact.
Together with Lemma~\ref{noline}, this leaves $M/G=[0,\infty).$

\begin{lemma}\label{monotone} Along the principal part
\[
 F<0,\qquad F'<0,\qquad H>0,\qquad H'<0
 \qquad\text{on }(0,\infty).
\]
\end{lemma}
\begin{proof}
By \eqref{criticalhess2}, one has $F<0$ and $ F'<0$ for $t>0$ sufficiently small. If $F$ were to vanish again, let $t_0>0$ be its first zero. Since $F<0$ on $(0,t_0)$, one would have $F'(t_0)\ge0$. On the other hand, the identity $F'+\Xi F=-C$
gives $F'(t_0)=-C<0,$ a contradiction. Hence $F<0$ on $(0,\infty).$

Suppose next that $F'$ has a first zero at $t_1>0$. At $t_1$, one has $R=-F'-HF=-HF>0.$ By $F<0$, it follows that $H(t_1)>0$, and therefore
\[
D(t_1)\ge\frac{H(t_1)^2}{3}>0.
\]
Using \eqref{radialid} and $F'(t_1)=0$, one obtains that
\[
F''(t_1)=D(t_1)F(t_1)<0.
\]
But at the first zero of a function $F'$ that is negative immediately beforehand one must have $F''(t_1)\ge0$, a contradiction. Thus $F'<0$ on $(0,\infty).$

It now follows from $H'=F'-D$ that $H'<0.$ Moreover, $H'\le -D\le-\frac{H^2}{3}.$ If $H$ ever became nonpositive, then, since $H'<0$, it would become strictly negative; the differential inequality above would force $H$ to diverge to $-\infty$ in finite forward time. This is incompatible with the existence of the complete principal normal geodesic on $[0,\infty)$. Consequently, $H>0$ on $(0,\infty).$
\end{proof}

Compare Lemma \ref{monotone} with \cite[Propositions~2.3 and~2.4]{BDW}, where these conclusions are obtained for complete cohomogeneity-one steady Ricci solitons with one special orbit.

\begin{lemma}
   The boundary orbit has codimension at least two.
\end{lemma}
\begin{proof}
We proceed by contradiction. If it has codimension one, then the slice isotropy acts by reflection in the normal direction, forcing the mean curvature of the nearby principal orbits to satisfy $H(0)=0.$ But $H'(0)=F'(0)-D(0)<0,$ so $H<0$ immediately outside the boundary, contradicting Lemma \ref{monotone}. 
\end{proof}
Next, since $F<0$ on $(0,\infty)$, there are no further critical orbits. Moreover, for every $t>t_0>0$,
\[
F(t)\le F(t_0)<0
\text{ and hence }
f(t)\longrightarrow-\infty
\text{ as }t\longrightarrow\infty.
\]

Finally, the principal part is a $G$-equivariant product over $(0,\infty)$. Attaching the slice at the singular orbit $N$ identifies $M$ with the normal bundle $\nu N$. More precisely, choose a unit-speed normal geodesic issuing from $N$. Its stabilizer is constant for $t>0$, it meets each principal orbit exactly once, and its $G$-translates yield the normal exponential map $\nu N\longrightarrow M.$ This map is a diffeomorphism: there can be no interior focal or identification points, since either phenomenon would produce an increase of isotropy or cause two distinct normal parameters to represent the same point of the orbit space. In particular, $N$ is connected and complete with its induced metric. Since $M$ deformation retracts onto the zero section of $\nu N$, and $M$ is simply connected, $N$ is simply connected as well.

\subsection{The point and surface alternatives}
We already know that the singular orbit $N$ has codimension $d\geq2$. We now distinguish the remaining possibilities.

\textbf{Case $d=3$:} Then $N$ is one-dimensional. Its connected isotropy fixes the tangent line pointwise, while acting transitively on the normal unit sphere $\mathbb S^2$. But it also fixes the image of that tangent line under $J$, and this image is a non-zero normal line, a contradiction. Thus $d=2$ or $4$.

\textbf{Case $d=4$:} Then $N$ is a fixed point. Properness makes $G$ compact, and its faithful isotropy representation is a four-dimensional connected closed subgroup of $U(T_NM,J)\cong U(2)$. It is therefore all of $U(2)$. The normal exponential identification yields $M\cong\R^4$, principal spheres $\mathbb S^3$, and the global biaxial Hopf metric, with the smooth nut data (\ref{nutdata}).

\textbf{Case $d=2$:} At $x\in N$, the isotropy has dimension two and acts on the normal two-plane by rotations. The identity component of the kernel of this normal action has dimension one. It acts non-trivially on $T_xN$ due to 
the effectiveness of the isometry action. Thus, its fixed subspace is exactly the normal two-plane. Since it commutes with $J$, both the normal and tangent two-planes are $J$-invariant. Thus $N$ is a complex curve, and its normal bundle is a complex line bundle.

The induced action of $G$ on $N$ has dimension three: its isotropy includes all rotations of the tangent two-plane, and a surface has isometry group dimension at most three. The kernel has a one-dimensional identity component $\mathbb S^1$. It is compact because it fixes $x$, and it is central because a connected group acts trivially by conjugation on the automorphism group of a circle. It rotates each normal fiber with weight one. Indeed, an element in the kernel of this rotation would fix both $x$ and its whole tangent space and hence be the identity isometry. This central circle fixes exactly $N$ and acts freely on $M\setminus N$. In particular, $N$ is totally geodesic.

Finally, homogeneity makes the Gaussian curvature of $N$ constant; completeness and simple connectivity make it a space form. 

\subsection{The bolt connection and positive transverse curvature}
Let $\mathcal P$ denote the unit normal circle bundle of $N$, and initially normalize the circle generator $K_0$ to have period $2\pi$. On each principal orbit, the orthogonal complement of the circle fibers defines a connection one-form $\eta$, normalized by $\eta(K_0)=1.$ The principal isotropy fixes the vertical line generated by $K_0$ and acts by rotations on the horizontal two-plane. Hence the invariant horizontal complement is unique. Under the $G$-equivariant normal-flow identification of the principal orbits, this uniqueness implies that $\eta$ is independent of $t$.

The same isotropy makes the horizontal metric a scalar multiple of a fixed metric $h$ on $N$. Its curvature is $d\eta=c\pi^*\omega_h$ for a constant $c$, since it is invariant under the full connected surface isometry group. We have therefore obtained the global form (\ref{ansatz}) on $(0,\infty)\times\mathcal P$, not just a local leaf-space model. At $t=0$,
\[
 p=t+O(t^3),\qquad q=q_0+O(t^2),\qquad q_0>0,
\]
with $p$ odd and $q,f,u$ even. 

\begin{lemma}\label{positivebase}
The transverse curvature $\kappa$ of a complete irreducible bolt solution is strictly positive, including in the untwisted case $c=0$.
\end{lemma}

\begin{proof} If $\kappa=c=0$, then (\ref{qweight}) gives $q'=0$. Here $N$ is the Euclidean plane, and the flat circle connection is globally trivial. Thus $g=(dt^2+p^2d\theta^2)+q_0^2h$ is a global Riemannian product, contrary to irreducibility.

Now assume $\kappa\leq0$ and $(\kappa,c)\neq(0,0)$. Two integrating-factor forms of (\ref{ode}) are
\begin{align}
 (q^2e^{-f}p')'&=\frac{c^2p^3}{2q^2}e^{-f},\label{pweight}\\
 (pqq'e^{-f})'&=p e^{-f}\left(\kappa-\frac{c^2p^2}{2q^2}\right).
 \label{qweight}
\end{align}
The first one gives $p'>0$. The second identity, with zero initial left-hand quantity, gives $q'<0$ for $t>0$. The scalar curvature of a principal orbit is
\[
 S_P(t)=\frac{2\kappa}{q^2}-\frac{c^2p^2}{2q^4}.
\]
Since $q\leq q_0$ and $p$ is increasing, there is a constant $\varepsilon>0$ such that $S_P\leq-\varepsilon$ for sufficiently large $t$. For $\kappa<0$, we use the first term; for $\kappa=0$ and $c\neq0$, we use $p(t)\geq p(t_0)>0$ in the second term. The scalar Gauss formula gives
\[
 R=S_P-H^2-D-2H'.
\]
Since $R>0$, $D\geq 0$, $H'\leq-\varepsilon/2$, contradicting $H>0$ on the entire half-line. 

\end{proof}
\begin{remark}
For twisted bolts, Lemma~\ref{positivebase} also follows from
\cite[Theorem~1.4]{TY}. That theorem asserts that $q\to\infty$ for every
complete steady GRS of this form. 
On the other hand, if
$\kappa\le0$, then \eqref{qweight} forces $q'<0$, so $q$ stays bounded.
\end{remark}
Consequently, $N$ is a round two-sphere. The kernel of the action of
$G$ on $N$ is compact, and its image is $SO(3)$, so $G$ is compact.
The normal bundle of $N$ is determined by its Euler number
$e\in\mathbb Z$. If the circle generator is normalized to have period
$2\pi$, then \[\int_N d\eta=c\,\operatorname{Area}(N,h)\] equals $2\pi e$
up to sign. Hence $e=0$ forces $c=0$. Then the connection is
flat; since $H^1(\mathbb S^2;\R)=0$, it is gauge equivalent to the product
connection $d\theta$ and \[M\cong \mathbb{S}^2 \times \mathbb R^2.\]

Suppose now that $m=|e|>0$. Choose the orientation of the fibers and
rescale $h$ so that $K_h=4$; then $\operatorname{Area}(N,h)=\pi$, and
the period-$2\pi$ connection has curvature $2m\,\omega_h$. Replacing
$\eta$ by $\eta/m$ gives $d\eta=2\omega_h$, fiber period $2\pi/m$ and
$p'(0)=m$. Two connections with the same curvature differ by a closed
one-form, which is exact because $H^1(\mathbb S^2;\R)=0$; hence they are gauge
equivalent. This is precisely the global lens-space ansatz with the bolt
normalization~\eqref{boltdata}. 


Furthermore, $G$ is a four-dimensional closed subgroup of the identity component of the group of isometries of $\mathcal P$ that preserve $K_0$. For the twisted metric $p^2\eta^2+q^2h$ on $\mathbb S^3/\mathbb Z_m$, that group is $U(2)/\mathbb Z_m$, acting through the standard action of $U(2)$ on $\C^2$. It is isomorphic to $U(2)$ for odd $m$ and to $SO(3)\times \mathbb S^1$ for even $m$. For the untwisted metric $p^2d\theta^2+q^2h$ on $\mathbb S^1\times \mathbb S^2$, it is $\mathbb S^1\times SO(3)$. Both groups are connected and four-dimensional, so $G$ is the whole group.
\subsection{Global one-variable data}

Since $f(t)$ is strictly decreasing on $(0,\infty)$, every regular level
of $f$ is a single connected $G$-orbit. As $G$ also preserves $u$, both
$|\nabla f|^2$ and $u$ are well-defined functions of $f$ on $(-\infty,0)$. 

It remains to check smoothness through the critical value $f=0$. The
slice isotropy acts transitively on the unit normal sphere of $N$, so
it contains an isometry reversing each normal geodesic issuing from
$N$. Hence every smooth $G$-invariant function is an even function of
the normal distance $t$, and by Whitney's theorem it is a smooth
function of $t^2$. We may therefore write
\[
 f(t)=\Phi(t^2),\qquad u(t)=\Psi(t^2),
\]
with $\Phi$ and $\Psi$ smooth on $[0,\infty)$. From~\eqref{criticalhess2}, one sees that
\[
 \Phi(0)=0,\qquad \Phi'(0)=\tfrac12 f''(0)=-\frac{C}{2d}<0 .
\]
Moreover, $\Phi'(s)=f'(t)/(2t)<0$ for $s=t^2>0$, so $\Phi$ is a
diffeomorphism from $[0,\infty)$ onto $(-\infty,0]$. Since
$|\nabla f|^2=f'(t)^2=4t^2\Phi'(t^2)^2$, the functions $B$ and $U$
defined by
\begin{equation}\label{Bendpoint}
 B(\Phi(s))=4s\,\Phi'(s)^2,\qquad U(\Phi(s))=\Psi(s)
\end{equation}
are smooth on $(-\infty,0]$, including the critical value, and satisfy
$|\nabla f|^2=B(f)$, $u=U(f)$ on $M$ and $B'(0)=4\Phi'(0)=-2C/d$. If
desired, they can be extended smoothly to all of $\R$. Finally, since
$\Delta f=-R=B(f)-C$, the transnormal data are also smooth and global.

This completes the proof of Proposition~\ref{globalsym}.

\section{Complete nut and bolt trajectories}\label{trajectories}

In this section, we classify the complete nut and bolt trajectories
of~\eqref{ode}. The uniqueness of the higher-degree threshold is
postponed to Section~\ref{thresholdsec}. Smoothness at the singular orbit means that $p$ and $q$ are odd
functions of $t$ at a nut, while at a bolt $p$ is odd and $q$ is
positive and even; in both cases $f$ is even. At a nut, we set
$\alpha=p^{(3)}(0)$ and $\beta=q^{(3)}(0)$. From this, it follows that
\begin{equation}\label{nutdata}
 \begin{gathered}
 p=t+\frac\alpha6t^3+O(t^5),\quad
 q=t+\frac\beta6t^3+O(t^5),\\
 f=f(0)+\frac{\alpha+2\beta}{2}t^2+O(t^4),\\
 C:=-4(\alpha+2\beta),\qquad
 Q:=\frac pq=1+\frac{\alpha-\beta}{6}t^2+O(t^4).
 \end{gathered}
\end{equation}
By~\cite[Lemma~3.9]{Donovan}, applied with $(f_1,f_2,f_3)=(p,q,q)$ and
$a_i=f_i^{(3)}(0)/3$, that is, $(a_1,a_2,a_3)=(\alpha/3,\beta/3,\beta/3)$,
these data determine a unique smooth germ. In particular, every smooth
biaxial nut germ is of this form.

At a bolt, we normalize $q(0)=1$ and $f(0)=0$, and set $\mu=f''(0)$.
Let $s_0=m$ for a twisted bolt of degree $m$ and $s_0=1$ for an
untwisted bolt. Then
\begin{equation}\label{boltdata}
 \begin{gathered}
 p=s_0t+\frac{s_0(\mu-4)}6t^3+O(t^5),\qquad
 q=1+t^2+O(t^4),\\ f=\frac\mu2t^2+O(t^4),\\
 C=-2\mu.
 \end{gathered}
\end{equation}

\begin{remark}\label{uniquegerm} The full germ is uniquely determined by $\mu$: see
\cite[Theorem~2.1]{Appleton} for twisted bolts and
\cite[Corollary~3.8, Theorem~4.7 and Remark~2.1]{Stolarski} for the
untwisted case; note that the potential in~\cite{Stolarski} is $-f$.
Since $C=R(0)$ and $R>0$, completeness and non-constancy force
$\alpha+2\beta<0$ at a nut and $\mu<0$ at a bolt.
\end{remark}

The next lemma proves the joint analytic dependence of twisted bolt germs, which will be crucial for our argument; compare \cite[Theorem~2.1 and Appendix~B]{Appleton}.  

\begin{lemma}\label{jointanalytic}
Fix $m\ge1$. For $\mu\in\R$, let $(p_\mu,q_\mu,f_\mu)$ be the twisted bolt solution
of~\eqref{ode} with data~\eqref{boltdata} of degree~$m$, and let $(0,T(\mu))$ be its maximal interval.
\begin{enumerate}[label=(\roman*)]
\item For every compact interval $I\subset\R$, there is $\varepsilon>0$
such that $T>\varepsilon$ on $I$ and
$(t,\mu)\mapsto\bigl(p_\mu(t),q_\mu(t),f_\mu(t)\bigr)$ extends to a
real-analytic map on a neighborhood of $[-\varepsilon,\varepsilon]\times I$;
here $p_\mu$ is odd and $q_\mu$, $f_\mu$ are even in $t$. Consequently,
\[
 p_\mu=mt+\frac{m(\mu-4)}6t^3+t^5\varphi_1,\qquad
 q_\mu=1+t^2+t^4\varphi_2,\qquad
 f_\mu=\frac\mu2t^2+t^4\varphi_3,
\]
where $\varphi_1,\varphi_2,\varphi_3$ are real analytic in $(t,\mu)$ on
that neighborhood. 
\item The function $T$ is lower semicontinuous, and
$(t,\mu)\mapsto\bigl(p_\mu(t),q_\mu(t),f_\mu(t)\bigr)$ is real analytic
on the open set of all $(t,\mu)\in\R^2$ with $0<t<T(\mu)$.
\end{enumerate}
\end{lemma}
\begin{proof}
(i) Since $p'>0$ near the bolt, we may use $r=p^2$ as the independent
variable, as in~\cite[Appendix~B]{Appleton}. Put $\mathfrak h=(p')^2$,
$y=dq/dr$ and $z=df/dr$, so that $q'=2pp'y$ and $f'=2pp'z$. For
$u=(\mathfrak h,q,y,z)$, the system~\eqref{ode} with $c=2$ and
$\kappa=4$ becomes
\[
 r\,\frac{du}{dr}=P(r,u),\qquad
 P(0,u)=\Bigl(0,\ 0,\ \frac1{\mathfrak hq}-y,\
 \frac2{\mathfrak hq^2}-\frac{2y}q\Bigr),
\]
where the components of $P$ are polynomials in
$r,\mathfrak h^{\pm1},q^{\pm1},y,z$. The data~\eqref{boltdata}
correspond to \[u(0)=u_0(\mu):=\bigl(m^2,1,m^{-2},\mu/(2m^2)\bigr), \text{ and }
P(0,u_0(\mu))=0.\] Since $P(0,u)$ does not involve $z$, the matrix
$A=\partial_uP(0,u_0(\mu))$ does not depend on $\mu$: its first two
rows vanish, and its last two rows are \[(-m^{-4},-m^{-2},-1,0) \text{ and }
(-2m^{-4},-2m^{-2},-2,0).\] Its eigenvalues are $0$ and $-1$, so
$kI-A$ is invertible for every integer $k\ge1$, and
$\sup_{k\ge1}\|(kI-A)^{-1}\|<\infty$ because $(kI-A)^{-1}\to0$.

The parameter-dependent Briot--Bouquet
theorem~\cite[Theorem~9.2]{Appleton}, applied to $v=u-u_0(\mu)$, gives
an analytic solution with $v(0,\mu)=0$, namely a convergent series
\[v=\sum_{k\ge1}v_k(\mu)r^k.\] Since $A$ is constant and $P$ is
polynomial in $z$, the recursion for the $v_k$ shows that they are
polynomials in $\mu$. Following the argument in~\cite{Appleton},
the series
converges uniformly on a complex neighborhood of $(0,\mu_0)$. By
compactness, $u$ is real analytic in $(r,\mu)$ on
$(-r_0,r_0)\times I'$ for some $r_0>0$ and some open interval
$I'\supset I$.

Since $\mathfrak h(0,\mu)=m^2>0$, the function
\[t(p,\mu)=\int_0^p\mathfrak h(s^2,\mu)^{-1/2}\,ds\] is real analytic and
odd in $p$, with $\partial_pt(0,\mu)=1/m$. By the analytic inverse
function theorem and the compactness of $I$, its inverse $p(t,\mu)$ is
real analytic on a neighborhood of $[-\varepsilon,\varepsilon]\times I$
for some $\varepsilon>0$, and odd in $t$. Together with $q(p^2,\mu)$
and $\int_0^{p^2}z(s,\mu)\,ds$, which are even in $t$, it solves
\eqref{ode} for small $t>0$ with the data~\eqref{boltdata}. By the
uniqueness of the germ (Remark~\ref{uniquegerm}), this solution is
$(p_\mu,q_\mu,f_\mu)$. 

After shrinking $\varepsilon$, we have
$p_\mu,q_\mu>0$ on $(0,\varepsilon]$, so $T>\varepsilon$ on $I$.
Finally, by the parities and~\eqref{boltdata}, the differences
$p_\mu-mt-m(\mu-4)t^3/6$, $q_\mu-1-t^2$, and $f_\mu-\mu t^2/2$
vanish to orders $5$, $4$, $4$ at $t=0$. Taylor's formula with integral
remainder in $t$ therefore gives the real-analytic functions
$\varphi_i$, whose derivatives are bounded on the compact set
$[-\varepsilon,\varepsilon]\times I$.

(ii) Let $\mu_0\in\R$, apply (i) with $I=[\mu_0-1,\mu_0+1]$, and fix
$t_0\in(0,\varepsilon)$. Then
$(p_\mu,p_\mu',q_\mu,q_\mu',f_\mu,f_\mu')(t_0)$ is real analytic in
$\mu\in I$, with $p_\mu(t_0),q_\mu(t_0)>0$. On $\{p>0,\ q>0\}$, the
first-order system associated with~\eqref{ode} is analytic, so by
analytic dependence on initial conditions its flow is real analytic on
an open domain. Hence $T$ is lower semicontinuous at $\mu_0$, and the
map in~(ii) is real analytic near every point of its domain with
$|\mu-\mu_0|<1$ and $t>t_0$. 
\end{proof}

\subsection{One continuation criterion}
Set $a=p'/p$ and $b=q'/q$, so that $H=a+2b$ and $D=a^2+2b^2$. By
\eqref{pweight}, the quantity $q^2e^{-f}p'$ is non-decreasing: it is
constant when $c=0$, and strictly increasing for $t>0$ when $c=2$. It
is positive at a bolt and vanishes at a nut. Hence $p'>0$ on every
maximal positive trajectory under consideration. By \eqref{boltdata}
and \eqref{nutdata}, $F<0$ and $F'<0$ for small $t>0$, since
$f''(0)=-C/d<0$. The identities \eqref{radialid} show that both signs
persist on the whole maximal positive interval. 

\begin{lemma}\label{continuation}
Let a maximal positive trajectory with $C>0$ satisfy
$q'>0$. Then it is complete. In particular, a twisted trajectory with
$Q<\sqrt2$ throughout is complete.
\end{lemma}

\begin{proof}
Write $(0,T)$ for the maximal positive interval. The preceding sign
argument gives $p'>0$, $F<0$ and $F'<0$. Since $q'>0$, we also have
$a,b,H>0$. Hence~\eqref{radialid} gives
\[
 R=-F'-HF>0,\qquad |F|<\sqrt C,\qquad
 H'=F'-D\le-\frac{H^2}{3}.
\]
At a nut $H\sim3/t$, and at a bolt $H\sim1/t$, so
$1/H\to0$ as $t\downarrow0$. Integrating $(1/H)'\ge1/3$ yields
\[
 0<H(t)\le\frac3t,\qquad 0<a(t),b(t)\le\frac3t
 \quad(0<t<T).
\]
These estimates hold on the maximal positive interval, before
completeness has been established.

Suppose $T<\infty$ and fix $t_0\in(0,T)$. Integrating the bounds on
$a=(\log p)'$ and $b=(\log q)'$ bounds $p,q$ above on $[t_0,T)$;
their monotonicity bounds them away from zero. The identities
$p'=ap$, $q'=bq$ and $|F|<\sqrt C$ also bound $p',q',f'$ there,
and integration bounds $f$. Thus the variables of the first-order
system associated to~\eqref{ode} stay in a compact subset of its
regular domain $\{p>0,q>0\}$. The solution extends past $T$, a
contradiction. Therefore $T=\infty$, and the smooth completion at
$t=0$, together with compact principal orbits, gives completeness.

For the last assertion, $q'>0$ initially by~\eqref{nutdata}
or~\eqref{boltdata}. A first zero of $q'$ while $Q<\sqrt2$ would
satisfy
\[
 q''=\frac{4-2Q^2}{q}>0,
\]
which is impossible after positive values of $q'$. Thus $q'>0$
throughout, and the first part applies; compare
\cite[Lemmas~4.1 and~5.3]{Appleton}.
\end{proof}

For twisted trajectories, the ratio equation is
\begin{equation}\label{ratio}
 Q''=(F-3b)Q'+\frac{4Q(Q^2-1)}{q^2}.
\end{equation}
A critical point in $0<Q<1$ is a strict maximum; one in $Q>1$ is a strict minimum. The next obstruction is needed for both endpoint types.

\begin{lemma}\label{Qobstruction}
Let $(p,q,f)$ be a complete non-trivial solution of~\eqref{ode} with
$c=2$ and $\kappa=4$, emanating from a nut or from a twisted bolt. Then
there is no $t_0>0$ with $Q(t_0)>1$ and $Q'(t_0)>0$.
\end{lemma}

For bolts, this is~\cite[Lemma~6.10]{Appleton}, combined
with~\cite[Remark~3.2]{Appleton}. We provide an argument that applies to both endpoint types.

\begin{proof}
The idea is that an increasing ratio above one must grow without bound.
This eventually makes the scalar curvature of the principal orbits
negative, which contradicts $H>0$ through the scalar Gauss formula.
Suppose, to the contrary, that $Q(t_0)>1$ and $Q'(t_0)>0$ for some
$t_0>0$.

By~\eqref{ratio}, at a zero of $Q'$ with $Q>1$ we have
$Q''=4Q(Q^2-1)/q^2>0$. Hence $Q'>0$ and $Q>1$ for all $t\ge t_0$. Also $p'>0$ by
\eqref{pweight}: the quantity $q^2e^{-f}p'$ vanishes at a nut and is positive at a bolt. 
Set
\[
 a=\frac{p'}{p}>0,\qquad b=\frac{q'}{q},\qquad F=f',\qquad w=qq'.
\]
The trajectory is complete, so the argument of
Subsection~\ref{subsection4.3} applies and gives 
\begin{align*}
    F<0,\, F'<0,\,|F| &\le\sqrt C \text{ and } H=a+2b>0 \text{ on } (0,\infty),\\
    H'=F'-D &\le-H^2/3 \text{ and } 1/H\to0 \text{ as } t\downarrow0,
\end{align*} at either type of singular orbit. Thus, $(1/H)'\ge1/3$ gives
\[
 0<H(t)\le\frac3t .
\]
The horizontal equation in~\eqref{ode} is equivalent to
\begin{equation}\label{wnut}
 w'+(a-F)w=4-2Q^2 .
\end{equation}
For $t\ge t_0$, set $\sigma=-F(t_0)>0$. Since $a>0$ and $F$ is
decreasing, $a-F\ge\sigma$. Let $M_{\max}=\max\{w(t_0),2/\sigma\}>0$. At any
time $t\ge t_0$ with $w\ge M_{\max}$, $Q>1$ gives
\[
 w'=4-2Q^2-(a-F)w<2-\sigma w\le0 .
\]
Hence $w\le M_{\max}$ on $[t_0,\infty)$. Since $(q^2)'=2w$, there is a
constant $K>0$ such that
\begin{equation}\label{qlinear}
 q(t)^2\le q(t_0)^2+2M_{\max}(t-t_0)\le K(1+t)\qquad(t\ge t_0).
\end{equation}
This bound does not use the sign of $q'$.

Since $a-b=Q'/Q>0$, we have $3b<a+2b=H\le3/t$. Therefore, for $t\ge t_0$,
\[
 -F+3b\le\sqrt C+\frac3{t_0}=:A_0 .
\]
Using $Q'>0$ in~\eqref{ratio}, the monotonicity of $Q$, and then
\eqref{qlinear}, we infer
\[
 Q''+A_0Q'\ \ge\ \frac{4Q(Q^2-1)}{q^2}\ \ge\ \frac{\gamma}{1+t},
 \qquad
 \gamma=\frac{4Q(t_0)\bigl(Q(t_0)^2-1\bigr)}{K}>0 .
\]
Thus $\bigl(e^{A_0t}Q'\bigr)'\ge\gamma e^{A_0t}/(1+t)$. Integrating over
$[t-1,t]$ using $Q'(t-1)>0$ gives
\[
 Q'(t)\ge\frac{\gamma e^{-A_0}}{1+t} \qquad(t\ge t_0+1).
\]
Consequently, $Q(t)\to\infty$ as $t\to\infty$.

Now, choose $T\ge t_0$ such that $Q\ge3$ on $[T,\infty)$. There,
\eqref{wnut} gives $w'\le4-2Q^2\le-14$ whenever $w\ge0$. Hence $w$
becomes negative at some time $T_1\ge T$ and stays negative afterwards,
since $w'<0$ at any zero of $w$. Therefore $q$ is decreasing on
$[T_1,\infty)$, and the scalar curvature of the principal orbits
satisfies
\[
 S_P=\frac{8-2Q^2}{q^2}\le-\frac{10}{q^2}\le-\frac{10}{q(T_1)^2}
 =:-\varepsilon<0, \qquad(t\ge T_1).
\]
By the scalar Gauss formula $R=S_P-H^2-D-2H'$ and $R>0$, we get
$H'<S_P/2\le-\varepsilon/2$ on $[T_1,\infty)$. This contradicts $H>0$
on the whole half-line. 
\end{proof}

Consequently, a twisted bolt is complete if and only if $Q<1$ throughout. Indeed, one direction is given by Lemma~\ref{continuation}. Conversely, suppose that $Q$
reaches $1$ for the first time at $t_1$, so that $Q'(t_1)\ge0$. If
$Q'(t_1)=0$, then $p=q$ and $p'=q'$ at $t_1$. The set $\{p=q\}$ is
invariant under~\eqref{ode} when $c=2$ and $\kappa=4$, so uniqueness
for the regular system would force $Q\equiv1$, contradicting $Q(0)=0$.
Hence $Q'(t_1)>0$, and Lemma \ref{Qobstruction} shows that the trajectory is
incomplete.

We are ready to classify nuts and low-degree bolts.

\begin{proposition}\label{completetraj}
Consider maximal trajectories of (\ref{ode}) with a smooth nut or bolt at $t=0$ and non-constant potential. Such a trajectory is determined by its germ data (\ref{nutdata}) or (\ref{boltdata}). If it is complete, then $p^{(3)}(0)+2q^{(3)}(0)<0$ at a nut and $\mu<0$ at a bolt. In the normalizations of Section~\ref{intro}, the following assertions hold.
\begin{enumerate}[label=(\roman*)]
\item A nut trajectory is complete if and only if $\delta\leq0$. The trajectory with $\delta=0$ is Bryant's soliton.
\item A twisted bolt trajectory with $\mu<0$ is complete if and only if $Q<1$ on its maximal positive interval.
\item For twisted bolts of degree one or two, and for untwisted bolts, every trajectory with $\mu<0$ is complete.
\end{enumerate}
\end{proposition}

\begin{proof}
We start with nuts. By~\eqref{nutdata}, $Q=1+\frac{\alpha-\beta}{6}t^2+O(t^4)$.
\begin{itemize}
\item If $\alpha>\beta$, then $Q>1$ and $Q'>0$ for small $t>0$, which is
excluded by Lemma~\ref{Qobstruction}.
\item If $\alpha<\beta$, then $Q<1$ and $Q'<0$ for small $t>0$. At a
first zero of $Q'$ with $0<Q<1$, \eqref{ratio} would give
$Q''=4Q(Q^2-1)/q^2<0$, which is impossible after negative values of
$Q'$. Hence $Q$ is decreasing and $0<Q<1$ throughout.
\item If $\alpha=\beta$, germ uniqueness gives $p\equiv q$, so $Q\equiv1$.
\end{itemize}
In the last two cases, Lemma~\ref{continuation} proves completeness.
After a homothety we may assume $\alpha+2\beta=-1$, that is, $C=4$, and
the remaining parameter is exactly $\delta=\alpha-\beta\le0$. The
endpoint $\delta=0$, where $p=q$, is Bryant's soliton~\cite{Bryant}.

For the natural complex structure, $g$ itself is K\"ahler exactly when
the conformal factor in~\eqref{conformal} is constant, that is, when
$qq'=p$. At a nut, the $t^3$-coefficients give $\frac23\beta=\frac16\alpha$,
that is, $\alpha=4\beta$. Conversely, $\alpha=4\beta$ gives a K\"ahler
germ by~\cite[Lemma~3.9]{Donovan}. With $\alpha+2\beta=-1$ this
gives $\beta=-1/6$ and $\delta=-1/2$, and by uniqueness this member is
H.-D.~Cao's soliton~\cite{Cao}. 

At a bolt, the following quantities are useful:
\[
 \zeta=\frac p{q^2},\qquad Z=qq'-p=q^2(b-\zeta),\qquad E=4-p(\Xi+a).
\]
A direct computation from~\eqref{ode} gives
\begin{equation}\label{ZEsystem}
 Z'=(F-a+2\zeta)Z+E,\qquad E'=2\zeta(b+\zeta)Z,
\end{equation}
and, at a bolt of degree $m$,
\[
 Z=(2-m)t+O(t^3),\qquad E(0)=4-2m .
\]

\emph{Degree two.} Put $Y=tZ$. Since $p$ is odd with $p'(0)=m$, the
function $a-1/t$ is smooth at $t=0$, and $\zeta=O(t)$. Thus
\eqref{ZEsystem} becomes a regular linear system for $(Y,E)$ with smooth coefficients:
\[
 Y'=\Bigl(F-\bigl(a-\tfrac1t\bigr)+2\zeta\Bigr)Y+tE,\qquad
 E'=\frac{2\zeta}{t}\,(b+\zeta)\,Y.
\]
For $m=2$ we have $Y(0)=0$ and $E(0)=0$, so
$Y\equiv E\equiv0$. Hence $qq'=p$, and every degree-two bolt is
K\"ahler for the natural complex structure. Moreover, $q'=p/q>0$, so
Lemma~\ref{continuation} shows that the trajectory is complete for
every $\mu<0$. These are the complete steady K\"ahler--Ricci solitons
on the canonical bundle of $\CP^1$ constructed by Dancer and
Wang~\cite[Theorem~4.20(i), with $n_1=0$ and $r=2$]{DW} and H.-D.~Cao~\cite{Cao}; see
also~\cite[Remark~6.6 and Theorem~4.7]{Stolarski}. Since $E(0)\ne0$
for $m\ne2$, no bolt of another degree is K\"ahler for the natural
complex structure.\\

\emph{Degree one.} Here $Z=t+O(t^3)>0$ and $E(0)=2$. As long as
$q'>0$, the two quantities stay positive together. At a first zero of
$Z$ we would have $Z'=E>0$, and $E'=2\zeta(b+\zeta)Z>0$ makes $E>2$ for
$t>0$. We claim that $Q<1/\sqrt2$ for all $t$. Suppose $Q$ first
reached $1/\sqrt2$ at $t_1$. By Lemma~\ref{continuation}, $q'>0$ on
$(0,t_1]$, since $Q<\sqrt2$ there. Moreover, $Q'(t_1)\ge0$ gives
$a\ge b$, and $Z>0$ gives $b>\zeta$. Since $F<0$,
\[
 p(\Xi+a)=p(2a+2b-F)>4p\zeta=4Q^2=2
\]
at $t_1$, that is, $E(t_1)<2$, a contradiction. Hence $Q<1/\sqrt2$
globally, and Lemma~\ref{continuation} proves completeness for every
$\mu<0$.\\ 

\emph{Untwisted bolts.} When $c=0$, \eqref{pweight} gives
\begin{equation}\label{untwisted}
 p'=\frac{e^f}{q^2}>0,
\end{equation}
because $q^2e^{-f}p'$ is constant and equals $1$ at $t=0$. A first zero
of $q'$ would satisfy $q''=4/q>0$, which is impossible. Hence $q'>0$,
and Lemma~\ref{continuation} proves completeness for every $\mu<0$. 
\end{proof}

\subsection{The higher-degree shooting partition}\label{shooting}
Fix $m\ge3$. For $\mu<0$, a \emph{shot} is the maximal positive
trajectory with the bolt data~\eqref{boltdata} of degree $m$, and
$Q_\mu$ denotes its ratio $Q$. We set
\begin{align*}
 \mathcal I_m&=\{\mu<0:\ Q_\mu>1\text{ at some finite regular time}\},\\
 \mathcal D_m&=\{\mu<0:\ 0<Q_\mu<1\text{ and }Q_\mu'<0
                \text{ at some finite regular time}\},\\
 \mathcal N_m&=(-\infty,0)\setminus(\mathcal I_m\cup\mathcal D_m).
\end{align*}
Both $\mathcal I_m$ and $\mathcal D_m$ are open, by continuous
dependence on $\mu$. They are disjoint. By~\eqref{ratio}, a critical
point of $Q$ in $(0,1)$ is a strict maximum, so once $Q'<0$ with
$0<Q<1$, the ratio decreases forever. After a crossing of one, $Q$
increases forever by the first part of the proof of Lemma~\ref{Qobstruction}.

Since $Q_\mu=mt+O(t^3)$, every shot starts with $0<Q<1$ and $Q'>0$. A
shot in $\mathcal N_m$ therefore keeps $0<Q<1$ and $Q'>0$ for all time:
it cannot reach $Q=1$, since a first arrival has $Q'>0$, and a zero of
$Q'$ in $(0,1)$ would be followed by $Q'<0$. A shot in $\mathcal D_m$
eventually has $Q'<0$. By Lemmas~\ref{continuation}
and~\ref{Qobstruction}, the set of complete shots is exactly
$\mathcal D_m\cup\mathcal N_m$.

We use three facts from~\cite{Appleton}:
\begin{itemize}
\item \cite[Theorem~5.1 and Remark~5.5]{Appleton} give a negative tail
of complete shots with $Q\le c_0<1$;
\item \cite[Lemma~7.1]{Appleton} gives a crossing of one at $\mu=0$, so
Lemma~\ref{jointanalytic} gives an interval of crossings next to
zero;
\item the bounded-fiber estimate~\cite[Lemma~6.9]{Appleton} when $Q\to0$ (Remark \ref{boundedfiber}).
\end{itemize}

Consider a complete trajectory, either a shot in
$\mathcal D_m\cup\mathcal N_m$ or one of the trajectories of Proposition \ref{completetraj}. The function $F$ is decreasing and satisfies
$-\sqrt C\le F<0$, so it has a limit $F_\infty<0$. Since
$0<H\le3/t$, the identity $F'=F^2-HF-C$ shows that $F'$ converges to
$F_\infty^2-C$. Because $F$ converges, this limit must be zero. Hence
\[
 F\longrightarrow-\sqrt C,\qquad R=C-F^2\longrightarrow0;
\]
compare~\cite[Corollary~2.6]{BDW} and~\cite[Corollary~3.7]{TY}.

Set $w=qq'$. The horizontal equation reads
\begin{equation}\label{weq}
 w'+(a-F)w=4-\frac{c^2}{2}Q^2 .
\end{equation}
On a complete trajectory, $q'>0$, so $w>0$; see the proof of
Lemma~\ref{continuation}. For $t\ge t_0$ we have $a-F\ge\sigma:=-F(t_0)>0$
and the right-hand side of~\eqref{weq} is at most $4$. Hence
$w\le\max\{w(t_0),4/\sigma\}$, and $q^2\le K(1+t)$ for some $K>0$.

Now consider a shot in $\mathcal D_m$, or a nut with $\delta<0$. After
some time $t_1$, its ratio is decreasing with $0<Q\le Q(t_1)<1$.
Suppose that $Q_\infty:=\lim Q>0$. Then $v=-Q'>0$ satisfies, by
\eqref{ratio},
\[
 v'+(-F+3b)v=\frac{4Q(1-Q^2)}{q^2}\ \ge\ \frac{k}{1+t} \qquad(t\ge t_1),
\]
with $k=4Q_\infty\bigl(1-Q(t_1)^2\bigr)/K>0$. The coefficient is bounded above: since $a>0$ and $H=a+2b\le 3/t$, we have
$2b<3/t$, so $$-F+3b\le A_0:=\sqrt C+9/(2t_1).$$ As $v>0$, this gives
$(e^{A_0t}v)'\ge ke^{A_0t}/(1+t)$, and integrating over $[t-1,t]$ yields
$v\ge ke^{-A_0}/(1+t)$ for $t\ge t_1+1$. Thus $\int^\infty v=\infty$,
contradicting $Q>0$, and so $Q\to0$.


For a shot in $\mathcal N_m$, $Q$ increases to a limit $Q_\infty\le1$.
If $Q_\infty<1$, then, since $F-3b<0$ and $Q'>0$,
\[
 Q''=(F-3b)Q'+\frac{4Q(Q^2-1)}{q^2}\le-\frac{k}{1+t}
\]
for large $t$ and some $k>0$. This forces $Q'\to-\infty$, which is
impossible. Therefore
\[
 Q\longrightarrow0\ \text{on }\mathcal D_m,\qquad
 Q\longrightarrow1\ \text{on }\mathcal N_m .
\]
In particular, the negative tail lies in $\mathcal D_m$, and an interval
$(\mu',0)$ lies in $\mathcal I_m$. Let $(-\infty,\mu_*)$ be the
connected component of $\mathcal D_m$ containing the tail. Then
$\mu_*\le\mu'<0$, and $\mu_*\notin\mathcal D_m$. Also
$\mu_*\notin\mathcal I_m$, because $\mathcal I_m$ is open and disjoint
from $\mathcal D_m$. Hence $\mu_*\in\mathcal N_m$, and $\mathcal N_m$
is non-empty.

Since $0<a\le3/t$, we have $a\to0$, and~\eqref{weq} with $a-F\to\sqrt C$
yields
\begin{equation}\label{endscales}
 qq'\longrightarrow\frac{4-2Q_\infty^2}{\sqrt C},\qquad
 \frac{q^2}{t}\longrightarrow\frac{8-4Q_\infty^2}{\sqrt C}
 \qquad(c=2).
\end{equation}
For degree one, $Z>0$ gives $0<p<qq'$, so $p$ is increasing and
bounded and converges to some $L\in(0,\infty)$. Moreover,
$qq'>p(t_0)>0$ for $t\ge t_0$, so $q\to\infty$ and $Q\to0$. Degree two
has $p=qq'$ and the same conclusions. For $c=0$, \eqref{untwisted}
gives $p'\le e^f$, which is integrable because $f\le F(t_0)(t-t_0)$,
and the right-hand side of~\eqref{weq} equals four. Thus every bolt of
degree one or two and every untwisted bolt has $p\to L\in(0,\infty)$
and $q^2\sim(8/\sqrt C)\,t$.

\begin{remark}\label{boundedfiber}
Equation~\eqref{endscales} recovers~\cite[Corollary~6.7 and
Lemma~6.8]{Appleton} in real dimension four. On every decreasing
twisted end, $Q_\infty=0$, so it gives
\[
 qq'\longrightarrow\frac4{\sqrt C},\qquad
 q^2\sim\frac8{\sqrt C}\,t,\qquad
 \frac{q'}q=\frac{qq'}{q^2}\sim\frac1{2t}.
\]
The bounded-fiber argument of~\cite[Lemma~6.9]{Appleton} uses these
estimates. As in the proof of~\cite[Theorem~8.1]{Appleton}, it therefore
gives $p\to L>0$ for all
decreasing twisted ends, including the nuts with $\delta<0$.
For the Bryant nut, $C=4$ and $p=q\sim\sqrt{2t}$.
\end{remark}
All ends are summarized in the following table.
\begin{equation}\label{ends}
 \begin{array}{c|c|c}
 \text{end}&p&q^2\\\hline
 \text{Bryant nut}&p=q&2t\,(1+o(1))\\
 \text{other nuts}&L+o(1)&4t\,(1+o(1))\\
 \text{collapsed bolts}&L+o(1)&(8/\sqrt C)\,t\,(1+o(1))\\
 \text{non-collapsed bolts}&p/q\to1&(4/\sqrt C)\,t\,(1+o(1))
 \end{array}
\end{equation}

The connection curvature formulas and~\eqref{ode} give
$|\operatorname{Rm}|=O(t^{-1})$:
\begin{itemize}
\item the orbit terms are controlled by $q^{-2}$ and $p^2/q^4$;
\item the mixed terms are controlled by $(p/q^2)|a-b|$;
\item the radial terms are controlled by $|a'+a^2|+|b'+b^2|$.
\end{itemize}
On an end with bounded fiber, a radial band of width $2\varepsilon\sqrt t$
around distance $t$ has volume $O(t^{3/2})$. For small fixed
$\varepsilon>0$, a ball of radius $r=\varepsilon\sqrt t$ centred in
that band satisfies $\sup|\operatorname{Rm}|\le r^{-2}$, and its volume
ratio is $O(t^{3/2})/r^4=O(t^{-1/2})\to0$. Hence these ends are
collapsed.

On an end with $p/q\to1$, write $\rho$ for the distance coordinate and
set $\sigma=(\rho-t)/\sqrt t$. After scaling the metric by $t^{-1}$,
the metrics are uniformly comparable on bounded $\sigma$-intervals to
\[
 d\sigma^2+\frac4{\sqrt C}\,(\eta^2+h)\qquad\text{on }\R\times L(m,1),
\]
so they have a uniform positive lower bound on volume ratios at these
and smaller scales. The horizontal sectional curvature
$(4-3Q^2)/q^2-b^2\sim\sqrt C/(4t)$ shows that every radius admissible
for the curvature bound is $O(\sqrt t)$, and compactness handles the
core. Hence these ends are non-collapsed.

\section{Uniqueness of the higher-degree threshold}
\label{thresholdsec}

Here we prove the complete range of twisted bolts
of degree $m\ge3$ in Theorem~\ref{classification}.

\begin{theorem}\label{threshold}
For each $m\geq3$, there is a unique $\mu_m<0$ such that
\[
 \mathcal D_m=(-\infty,\mu_m),\qquad
 \mathcal N_m=\{\mu_m\},\qquad
 \mathcal I_m=(\mu_m,0).
\]
In particular, the unique non-collapsed normalized bolt is the endpoint of the complete interval $(-\infty,\mu_m]$.
\end{theorem}

The strategy is as follows. Shots in $\mathcal N_m$ stay in a simple barrier region $\mathcal R$ (Lemma ~\ref{barrier}). On $\mathcal R$, the derivatives of a shot with respect to $\mu$ satisfy a linear system whose off-diagonal coefficients are non-negative. They start positive by (\ref{variationjets}) and stay positive (Lemma~\ref{variational}). Two thresholds would then make $\log(Q_{\mu_2}/Q_{\mu_1})$ positive and increasing, which is incompatible with $Q\to1$ along both shots. 

Recall that $F=f'$, $a=p'/p$, $b=q'/q$, $H=a+2b$ and
$\Xi=H-F$. Set
\[
 \vartheta=a-b=\frac{Q'}{Q},\qquad
 \chi=\frac{qQ'}{1-Q^2},\qquad
 \mathcal R=\{0<Q<1,\ \chi>1\}.
\]
The ratio equation~\eqref{ratio} gives
\begin{equation}\label{chieq}
 \chi'=(F-2b)\chi+\frac{2Q}{q}\,(\chi^2-2),\qquad \chi(0)=m,
\end{equation}
where $\chi(0)=m$ follows from $Q=mt+O(t^3)$ and $q(0)=1$.

\begin{lemma}\label{barrier}
Every shot in $\mathcal N_m$ stays in $\mathcal R$.    
\end{lemma}
\begin{proof}
Indeed, suppose it reached $\chi=1$ at $t_0$. Then $Q'> 0$ and \eqref{chieq} imply $0<\chi\leq1$ thereafter, and $\chi'\leq F(t_0)\chi$. But $q\geq1$, so
\[
 (\operatorname{arctanh}Q)'=\chi/q\leq\chi
\]
would be integrable, contradicting $Q\to1$.
\end{proof}
The background equations are
\begin{equation}\label{background}
 \begin{gathered}
 a'=-\Xi a+2p^2/q^4,\qquad
 b'=-\Xi b+4/q^2-2p^2/q^4,\qquad \Xi'=-a^2-2b^2,\\
 H'+\Xi H=S_P,\qquad S_P=(8-2Q^2)/q^2,\\
 \vartheta'=-\Xi\vartheta-4(1-Q^2)/q^2.
 \end{gathered}
\end{equation}
On $\mathcal R$, one has $H>0$, $H'<0$, $\vartheta>0$, and $\vartheta'<0$.

Next, we differentiate at fixed $t$ with respect to $\mu$, keeping $m$ and $q(0)=1$ fixed, and define
\begin{equation}\label{variations}
 \begin{gathered}
 x_p=\partial_\mu\log p,\quad x_q=\partial_\mu\log q,\quad
 \mathcal V=x_p+2x_q,\quad \mathcal Z=x_p-x_q,\quad
 \mathcal U=-\partial_\mu\Xi,\\
 \mathcal T=\mathcal V'+\tfrac23H\mathcal V,\quad
 \mathcal W=\mathcal Z',\quad
 \mathcal J=\mathcal U-\tfrac23H\mathcal V-\tfrac43\vartheta\mathcal Z,\\
 \gamma=\frac{8(1-Q^2)}{3q^2},\qquad k_0=\frac{8(4Q^2-1)}{3q^2}.
 \end{gathered}
\end{equation}
Differentiating~\eqref{background} with respect to $\mu$ therefore gives
\[
 \begin{aligned}
 \mathcal V''+\Xi\mathcal V'+\tfrac23S_P\mathcal V
      &=H\mathcal U+2\gamma\mathcal Z,\\
 \mathcal Z''+\Xi\mathcal Z'
      &=\vartheta\mathcal U+\gamma\mathcal V+k_0\mathcal Z,\\
 \mathcal U'&=\tfrac23H\mathcal V'+\tfrac43\vartheta\mathcal Z'.
 \end{aligned}
\]
Substituting (\ref{variations}), and using $H'+\Xi H=S_P$, yields

\begin{equation}\label{firstorder}
 \begin{aligned}
 \mathcal V'&=-\tfrac23H\mathcal V+\mathcal T,\\
 \mathcal T'&=-(\Xi-\tfrac23H)\mathcal T+\tfrac29H^2\mathcal V
       +(2\gamma+\tfrac43H\vartheta)\mathcal Z+H\mathcal J,\\
 \mathcal Z'&=\mathcal W,\\
 \mathcal W'&=-\Xi\mathcal W+(k_0+\tfrac43\vartheta^2)\mathcal Z
       +(\gamma+\tfrac23H\vartheta)\mathcal V+\vartheta\mathcal J,\\
 \mathcal J'&=-\tfrac23H'\mathcal V-\tfrac43\vartheta'\mathcal Z.
 \end{aligned}
\end{equation}
\begin{lemma}\label{variational}
  On $\mathcal R$, every off-diagonal coefficient is non-negative.  
\end{lemma}
 \begin{proof}
For all but one of them, this follows from the signs of $H$, $H'$, $\vartheta$, and $\vartheta'$ on $\mathcal R$ recorded above, together with $\gamma>0$. Only the coefficient $k_0+\frac43\vartheta^2$ needs checking: for $s=Q^2\in(0,1)$ and $\chi>1$,
\[
 \begin{aligned}
 q^2(k_0+\tfrac43\vartheta^2)
 &=\tfrac43\left[2(4s-1)+\chi^2(1-s)^2/s\right]\\
 &>\tfrac43\left[(3s-1)^2/s+2\right]>0.
 \end{aligned}
\]
\end{proof}

Next, we determine the initial signs. The bolt expansion (\ref{boltdata}) gives
\begin{equation}\label{variationjets}
 \begin{gathered}
 \mathcal V=\tfrac16t^2+O(t^4),\qquad
 \mathcal Z=\tfrac16t^2+O(t^4),\qquad
 \mathcal W=\tfrac13t+O(t^3),\\
 \mathcal T=\tfrac49t+O(t^3),\qquad
 \mathcal J=\tfrac13t+O(t^3).
 \end{gathered}
\end{equation}
For instance, $\Xi=t^{-1}+(8-2\mu)t/3+O(t^3)$, so $\mathcal U=2t/3+O(t^3)$. 

By the joint
analyticity in $(t,\mu)$, the remainders and their $\mu$-derivatives
are uniform on compact parameter intervals. Hence there is $t_*>0$,
uniform on such intervals, with all five variables positive on
$(0,t_*]$.

The system~\eqref{firstorder} has continuous coefficients on every
compact interval in $(0,\infty)$ on which the shot remains in
$\mathcal R$. Write it as $y_i'=c_{ii}y_i+\sum_{j\ne i}c_{ij}y_j$ with
$c_{ij}\ge0$ for $j\ne i$. If all $y_j\ge0$ up to some time, then
$\bigl(e^{-\int c_{ii}}y_i\bigr)'\ge0$, so no component can reach zero
first. Starting at $t_*$, it follows that, as long as the shot remains
in $\mathcal R$,
\[
 \mathcal V,\mathcal T,\mathcal Z,\mathcal W,\mathcal J>0,
 \qquad \partial_\mu Q>0,
 \qquad \partial_t\partial_\mu\log Q>0.
\]
Here we use $\partial_\mu\log Q=\mathcal Z$ and $\partial_t\partial_\mu\log Q=\mathcal W$, by (\ref{variations}). Finally,
\[
 \partial_\mu\log\chi
 =\tfrac13\mathcal V+
   \frac{2(1+2Q^2)}{3(1-Q^2)}\mathcal Z+
   \frac{\mathcal W}{\vartheta}>0.
\]

We are ready to prove Theorem~\ref{threshold}.

\begin{proof}[{\bf Proof of Theorem~\ref{threshold}}]
The idea is to compare two shots in $\mathcal N_m$ through the whole
family of shots between them. Monotonicity in the parameter makes the
quotient of their ratios $Q$ increase, while both ratios tend to one.

Suppose that $\mu_1<\mu_2$ both belong to $\mathcal N_m$. By
Lemma~\ref{barrier}, their shots stay in $\mathcal R$ for all $t>0$.
By Lemma~\ref{jointanalytic}, there
is $t_*>0$ such that every shot with $\mu\in[\mu_1,\mu_2]$ lies in
$\mathcal R$ on $(0,t_*]$. Let $T_*$ be the supremum of all $T>0$
such that every shot with $\mu\in[\mu_1,\mu_2]$ exists and lies in
$\mathcal R$ on $(0,T]$. 

For $0<t<T_*$, Lemma~\ref{variational}
gives $\partial_\mu Q>0$ and $\partial_\mu\log\chi>0$ along the whole
family. Integrating in $\mu$, we obtain
\begin{align*}
 0 &<Q_{\mu_1}(t)\le Q_\mu(t)\le Q_{\mu_2}(t)<1,    \\
 1 &<\chi_{\mu_1}(t)\le\chi_\mu(t)\le\chi_{\mu_2}(t) \qquad(\mu_1\le\mu\le\mu_2). 
\end{align*}

Suppose that $T_*<\infty$. We first show that the whole family
continues past $T_*$.
\begin{itemize}
\item Since $C=-2\mu\le-2\mu_1$, we have $|F|\le\sqrt{-2\mu_1}$ for all
$\mu\in[\mu_1,\mu_2]$.
\item On $\mathcal R$ we have $q'>0$, so
$0<a,\,b\le3/t$.
\item Hence, on $[t_*,T_*)$, the functions $p$ and $q$ and their first
derivatives are bounded, and $p$ and $q$ are bounded away from zero,
uniformly in $\mu$.
\end{itemize}
The first-order system associated to~\eqref{ode} is regular in this
range, so every shot extends beyond $T_*$, with joint continuity in
$(t,\mu)$ by Lemma~\ref{jointanalytic}. At $t=T_*$, the inequalities above and
Lemma~\ref{barrier} give
\[
 Q_\mu(T_*)\le Q_{\mu_2}(T_*)<1,\qquad
 \chi_\mu(T_*)\ge\chi_{\mu_1}(T_*)>1
 \qquad(\mu_1\le\mu\le\mu_2),
\]
so the shots satisfy these strict inequalities uniformly in $\mu$. By
joint continuity and the compactness of $[\mu_1,\mu_2]$, all shots
remain in $\mathcal R$ on $(0,T_*+\varepsilon]$ for some
$\varepsilon>0$. This contradicts the definition of $T_*$, so
$T_*=\infty$.

Consequently, for every $t>0$,
\[
 G(t):=\log\frac{Q_{\mu_2}(t)}{Q_{\mu_1}(t)}
     =\int_{\mu_1}^{\mu_2}\mathcal Z_\mu(t)\,d\mu>0,\qquad
 G'(t)=\int_{\mu_1}^{\mu_2}\mathcal W_\mu(t)\,d\mu>0,
\]
where differentiation under the integral sign is justified by the joint
analyticity in $(t,\mu)$. Thus $G$ is positive and increasing, so
$G(t)\ge G(1)>0$ for $t\ge1$. On the other hand, both $Q_{\mu_1}$ and
$Q_{\mu_2}$ tend to $1$, by Subsection~\ref{shooting}, so $G(t)\to0$.
This is a contradiction.

Therefore $\mathcal N_m$ contains at most one point. It is non-empty by
Subsection~\ref{shooting}; write $\mathcal N_m=\{\mu_m\}$. The interval
$(-\infty,\mu_m)$ is connected and is the disjoint union of the open
sets $\mathcal D_m\cap(-\infty,\mu_m)$ and
$\mathcal I_m\cap(-\infty,\mu_m)$. It contains the negative tail, which
lies in $\mathcal D_m$, so $(-\infty,\mu_m)\subset\mathcal D_m$.
Likewise, $(\mu_m,0)$ contains an interval next to zero, which lies in
$\mathcal I_m$, so $(\mu_m,0)\subset\mathcal I_m$. This proves the
asserted partition.

Finally, the complete bolts of degree $m$ are those with
$\mu\in\mathcal D_m\cup\mathcal N_m=(-\infty,\mu_m]$. By~\eqref{ends},
the bolts with $\mu<\mu_m$ are collapsed and the bolt with $\mu=\mu_m$
is non-collapsed. This concludes the proof of Theorem~\ref{threshold}.
\end{proof}

\section{The converse and intrinsic parameters}\label{converse}

In this section, we verify the converse part of
Theorem~\ref{classification} and separate the parameters intrinsically.
We then complete the proof of Theorem~\ref{classification}. 
\subsection{The Converse} \label{converse1}
Let $K$ be the generator of the fiber circle, so that $\eta(K)=1$.
Each completed model carries a natural almost complex structure $J$: it sets $J\partial_t=p^{-1}K$ and acts on the horizontal distribution
as the lift of the complex structure of $(N,h)$. Its fundamental form
is
\[
 \omega_g=p\,dt\wedge\eta+q^2\omega_h .
\]
The $(1,0)$-forms of $J$ are spanned by $dt+ip\,\eta$ and the pull-back
of a $(1,0)$-form of the base. Since
$d(dt+ip\,\eta)=i\bigl(p'\,dt\wedge\eta+cp\,\omega_h\bigr)$, and both
$dt\wedge\eta$ and $\omega_h$ have type $(1,1)$, these derivatives have
no $(0,2)$-part. Hence $J$ is integrable. Moreover,
$$d(e^{2u}\omega_g)=\bigl((e^{2u}q^2)'-ce^{2u}p\bigr)\,dt\wedge\omega_h,$$
so $e^{2u}g$ is K\"ahler exactly when
\begin{equation}\label{conformal}
 (e^{2u}q^2)'=ce^{2u}p,
 \qquad\text{that is,}\qquad
 u'=\frac{cp}{2q^2}-\frac{q'}{q}.
\end{equation}

We check that $u$ is smooth through the singular orbit.
\begin{itemize}
\item \emph{At a nut}, write $p=tP(t^2)$ and $q=tQ_1(t^2)$, with
$P(0)=Q_1(0)=1$. Then
\[
 \frac{p}{q^2}-\frac{q'}{q}
 =\frac1t\Bigl(\frac{P}{Q_1^2}-\frac{Q_1+2t^2Q_1'}{Q_1}\Bigr).
\]
The bracket is smooth, even and vanishes at $t=0$, so the right-hand
side is a smooth odd function of $t$.
\item \emph{At a bolt}, $p$ is odd and $q$ is positive and even, so
both terms of $u'$ are smooth and odd.
\end{itemize}
In both cases $u$ is a smooth even function of $t$, hence a smooth
function on $M$, and $e^u>0$. For $c=0$ we have
$u=-\log q+\mathrm{const}$, and
$e^{2u}g=q^{-2}(dt^2+p^2\eta^2)+h$ is a product of two surfaces. 

The complex structure extends smoothly over the singular orbit. Near a
twisted bolt, put $d\log\rho=m\,dt/p$; near a nut or an untwisted bolt,
put $d\log\rho=dt/p$. In both cases $\rho\sim t$, and $\rho e^{i\psi}$
is a holomorphic fiber coordinate, where $\psi$ is the fiber angle
normalized to period $2\pi$.

Moreover, the potential satisfies $f'<0$ for $t>0$ and $f''(0)<0$, and
it has even expansions at the singular orbit. The inverse-function
argument of~\eqref{Bendpoint} therefore shows that $|\nabla f|^2=B(f)$
and $\Delta f=B(f)-C$ are smooth global functions of $f$. Hence $f$ is
leafwise isoparametric.

All constructed manifolds, namely $\R^4$, $\OO(-m)$ and
$\mathbb{S}^2 \times \mathbb R^2$, are simply connected, and all constructed metrics are
complete. By de Rham's theorem, a reducible metric would split off
either a flat factor, which yields a non-zero parallel vector field, or
a product of two non-flat surfaces.

\emph{Nuts.} At a normalized nut, $\Ric=-\nabla^2f=g$, so no parallel
vector field exists: it would lie in the kernel of $\Ric$. A product of
two non-flat surfaces would have a splitting of $T_0M$ into two
two-planes, which the connected group $U(2)$ would preserve, a contradiction. 

\emph{Bolts.} The normalized bolt $N$ is totally geodesic, and at $N$
\[
 \begin{gathered}
 K_N=4,\qquad \Ric|_{\nu N}=-\mu\,g,\qquad \Ric|_{TN}=0,\\
 K(\partial_t,X)=-\frac{q''(0)}{q(0)}=-2\qquad(X\in TN,\ |X|=1).
 \end{gathered}
\]
A parallel vector field $V$ would satisfy $\Ric(V)=0$, so $V\in TN$ at
$N$, and $R(\cdot,\cdot)V=0$. In particular, the sectional curvature of
the tangent plane $TN=\operatorname{span}\{V,JV\}$ would vanish,
contradicting $K_N=4$. In a product of two surfaces, the factor planes
at a point are Ricci eigenspaces. Since the eigenvalues $0$ and
$-\mu>0$ are distinct, the factor planes at $N$ would be $TN$ and
$\nu N$, and mixed sectional curvatures would vanish. This contradicts
$K(\partial_t,X)=-2$. 

\subsection{Separation of the parameters}\label{separation}
On an irreducible soliton, the potential is unique up to an additive
constant. Therefore, an isometry between two solitons maps the critical orbit of one potential onto that of
the other. The topology of the critical orbit and the value of
the Euler number of its normal bundle distinguish nuts from bolts,
twisted from untwisted bolts, and bolts of different degrees.

\emph{Nuts.} The condition $R(0)=C=4$ fixes the normalized scale. For
every unit vector $\nu\in T_0M$, the Jacobi operator $R(\cdot,\nu)\nu$
on $\nu^\perp$ has eigenvalues
\[
 -\alpha=\frac{1-2\delta}{3}\quad\text{(multiplicity one)},\qquad
 -\beta=\frac{1+\delta}{3}\quad\text{(multiplicity two)}.
\]
These are the limits of $-p''/p$ and $-q''/q$ as $t\to0$, and $U(2)$
acts transitively on unit directions. For $\delta<0$, the simple
eigenvalue determines $\delta$, while for $\delta=0$ the three
eigenvalues coincide. Hence distinct normalized nuts are not
homothetic.

\emph{Bolts.} The curvature $K_N=4$ of the critical orbit fixes the
scale, and the scale-invariant quantity
\[
 \frac{R|_N}{K_N}=-\frac\mu2
\]
distinguishes the normalized parameters. 

We are ready to prove Theorem~\ref{classification}.

\begin{proof}[{\bf Proof of Theorem~\ref{classification}}]
The proof assembles the previous results. If $W_g\not\equiv0$, Theorem~\ref{alignment}
with $\lambda=0$ gives $du\wedge df\equiv0$. If $W_g\equiv0$, then
$g$ is not flat, by irreducibility, so it is the Bryant soliton ~\cite[Theorem~1.2]{CC}; see also~\cite{CM}. In that case
Proposition~\ref{bryant} gives $du\wedge df\equiv0$.

Then, by Proposition~\ref{normalform}, every
regular point has a neighborhood with a four-dimensional Lie algebra of
Killing fields preserving $f$, $u$ and $J$. By
Proposition~\ref{globalsym}, which uses completeness and simple
connectivity, $(M,g,f)$ is a global nut or bolt trajectory of the form
\eqref{ansatz}, with the endpoint and period conventions of
Section~\ref{intro}.

Finally, Proposition~\ref{completetraj} and Theorem~\ref{threshold} give exactly the
parameter ranges in the statement. The converse follows from Subsection \ref{converse1}. 
By Subsection~\ref{separation}, distinct normalized parameters and distinct rows or degrees give
non-homothetic metrics. The nut $\delta=0$ is Bryant's soliton by Proposition~\ref{completetraj}, and the statements on collapse follow
from~\eqref{ends}, the discussion after it, and
Theorem~\ref{threshold}.
\end{proof}

\bibliographystyle{plain}
\bibliography{refs}

\end{document}